\documentclass{article}
\usepackage[latin1]{inputenc}
\usepackage{graphicx}
\usepackage{pgf,tikz}
\usetikzlibrary{arrows}
\usepackage{epstopdf}
\usepackage{dsfont}

\usepackage{amsfonts}
\usepackage{amsmath, latexsym}
\usepackage{amssymb,verbatim}
\usepackage{theorem}
\usepackage{color}
\usepackage{euscript}

\newtheorem{thm}     {Theorem}
\newtheorem{lem} [thm] {Lemma}
\newtheorem{pro}  [thm]   {Proposition}

\newtheorem{cor} [thm]    {Corollary}
\newtheorem{df}   [thm]   {Definition}

\theorembodyfont{\rmfamily}
\newtheorem{rem}[thm] {Remark}

\newcommand{\us}{u^\star}
\newcommand{\ps}{\psi^\star}

\newcommand{\om}{\omega}
\newcommand{\ba}{\begin{align*}}
\newcommand{\ea}{\end{align*}}

\newcommand{\Se}{{\cal S}}

\newcommand{\RR}{\mathbb{R}}
\newcommand{\R}{\mathbb{R}}
\newcommand{\Ss}{\mathbb{S}}

\newcommand{\eps}{\varepsilon}
\newcommand{\rti}{\tilde{r}}
\newcommand{\tit}{\tilde{t}}

\newcommand{\h}{{\mc{H}}}

\newcommand{\proof}[1]{\par\medskip\noindent{\bf Proof#1.}}
\newcommand{\qed}{\hfill$\square$}

\newcommand{\be}{\begin{equation}}

\newcommand{\ee}{\end{equation}}

\newcommand{\loc}{_{loc}}

\newcommand{\nd}{\noindent}

\newcommand{\Om}{\Omega}

\newcommand{\abs}[1]{\left\lvert#1\right\rvert}

\newcommand{\f}{\varphi}

\newcommand{\supp}{\operatorname{supp}}

\newcommand{\diff}{\, d}

\def\I[#1]{\mc{I}_{#1}}

\def\E[#1]{\mc{E}_{#1}}
\def\S[#1]{\mc S_0(#1)}
\def\mc{\mathcal}

\def\ed0{\eps\downarrow0}

\title{Kinetic formulation of vortex vector fields
}

\author{ {\Large Pierre Bochard} \footnote{D\'epartement de Math\'ematiques, Universit\'e Paris-Sud 11, 91405 Orsay, France. Email: Pierre.Bochard@math.u-psud.fr} \and
{\Large Radu Ignat}
\footnote{Institut de Math\'ematiques de Toulouse, Universit\'e
Paul Sabatier, 31062
Toulouse, France. Email: Radu.Ignat@math.univ-toulouse.fr}
}

\begin{document}

\maketitle

\begin{abstract}
This article focuses on gradient vector fields of unit Euclidean norm in $\R^N$. The stream functions 
associated to such vector fields solve the eikonal equation and the prototype is given by the distance function to a closed set.
We introduce a kinetic formulation that characterizes stream functions whose level sets are either spheres or hyperplanes in dimension $N\geq 3$. Our main result proves that the kinetic formulation is a selection principle for the vortex vector field whose stream function is the distance function to a point. 
\end{abstract}


\section{Introduction}
\label{chap_div}

In this article, we analyze the following type of vortex vector field:
$$\us:\RR^N\to \RR^N, \quad \us(x)=\frac{x}{|x|} \quad \textrm{ for every } \quad x\in \R^N\setminus\{ 0\}$$ 
in dimension $N\geq 2$ where $|\cdot|$ is the Euclidian norm in $\RR^N$. This structure arises in many physical models such as micromagnetics, liquid crystals, 
superconductivity, elasticity. Clearly, $\us$ is smooth away from the origin: in fact, $0$ is a topological singularity
of degree one since the jacobian is $\det \nabla \us=V_N {\bf \delta}_0$ where ${\bf \delta}_0$ is the Dirac measure at the origin and $V_N$ is the volume of the unit ball in $\RR^N$. Also, $\us$ is a curl-free 
unit-length vector field, i.e., 
\be
\label{eq:contraint}
|\us|=1  \quad \textrm{and} \quad \nabla \times \us=0 \quad \textrm{in } \, \RR^N\setminus \{0\}.
\ee
Moreover, there is a stream function $\ps:\RR^N\to \RR$ associated to $\us$ by the equation $$\us=\nabla \ps;$$ indeed, one may consider $\ps$ as the distance 
function at the origin, i.e., $\ps(x)=|x|$ for $x\in \RR^N$ and $\ps$ represents the viscosity solution of the eikonal equation
$$|\nabla \ps|=1$$
under an appropriate boundary condition at infinity (e.g., $\lim_{|x|\to \infty} ({\ps(x)}-{|x|})=0$).

Note that conversely, these properties characterize the vortex vector field: if $u:\RR^N \to \RR^N$ is a nonconstant vector field that is smooth away 
from the origin and satisfies \eqref{eq:contraint} then
$u=\pm \us$ in $\RR^N$. Indeed, this classically follows by the method of characteristics: the flow associated to $u$ by
\be
\label{eq:syst_dyn}
\partial_t{X}(t,x)=u(X(t,x))
\ee with the initial condition $X(0,x)=x$ for $x\neq 0$ yields straight lines $\{X(t, x)\}_t$ given by $X(t,x)=x+t u(x)$ along which $u$ is constant, i.e., 
$u(X(t,x))=u(x)$. Since $u$ is nonconstant and two characteristics can intersect only at the origin (which is the prescribed point-singularity of $u$), then every characteristic passes through the origin \footnote{This argument is clear in dimension $N=2$; for dimensions $N\geq 3$, one needs an additional argument showing that two characteristics are coplanar as we will see later in the proof of Theorem \ref{thm:main}.}
and therefore, $u$ coincides with $\us$ or $-\us$. In a recent paper, Caffarelli-Crandall \cite{CafCra} proved this result under a weaker regularity hypothesis for the vector field $u=\nabla \psi$: if $\psi$ is assumed only pointwise differentiable away from a set $S$ of vanishing Hausdorff $\h^1$-measure (i.e., $\h^1(S)=0$) and $|\nabla \psi|=1$ in $\RR^N\setminus S$, then $\psi=\pm \ps$ (up to a translation and an additive constant). We also refer to DiPerna-Lions~\cite{DipLio} for weaker regularity assumptions on $u$ in the framework of Sobolev spaces.

Our aim is to prove a kinetic characterization of the vortex vector field that does not assume any initial regularity on $u$. 
This kinetic formulation will characterize stream functions whose level sets are totally umbilical hypersurfaces in dimension $N\geq 3$, i.e., either pieces of spheres or hyperplanes. In order to introduce the kinetic formulation of the vortex vector field, we start by presenting the case of dimension $N=2$ and then we extend it to dimensions $N\geq 3$.

\subsection{Kinetic formulation in dimension $N=2$}

\nd Let $\Omega\subset \RR^2$ be an open set and $u:\Omega\to \RR^2$ be a Lebesgue measurable vector field that satisfies
\be
\label{eq:contraint2}
|u|=1 \, \textrm{ a.e. in } \Omega \quad \textrm{ and } \quad \nabla \times u=0 \, \textrm{ distributionally in } \Omega.
\ee
The main feature of the kinetic formulation relies on the concept of weak characteristic for a nonsmooth vector field $u$. 
\begin{figure}[htbp]
\center
\includegraphics[scale=0.3,
width=0.3\textwidth]{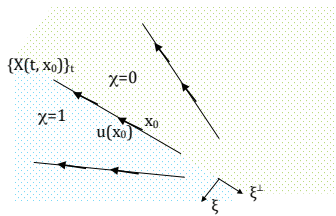} \caption{Characteristics of $u$.} \label{chara}
\end{figure}
We start by noting that \eqref{eq:syst_dyn} has a proper meaning only if some notion of trace of $u$ can be defined on curves $\{X(t, x)\}_t$ 
which in general is a consequence of the regularity assumption on $u$ (see DiPerna-Lions \cite{DipLio}). To overcome this difficulty, the 
following notion of ``weak characteristic" is introduced for measurable vector fields $u$ (see e.g. Lions-Perthame-Tadmor \cite{LPT}, Jabin-Perthame
\cite{Jabin-Perthame}): for every direction $\xi\in \Ss^1$, one defines the function 
$\chi(\cdot, \xi):\Omega\to\{0, 1\}$ by
\be
\label{eq:weak_cara}
\chi(x, \xi)=\begin{cases} 1 &\quad \textrm{ for }\, \, u(x)\cdot \xi>0, \\
0 &\quad \textrm{ for }\, \, u(x)\cdot \xi\leq 0.
\end{cases}\ee
In the case of a smooth vector field $u$ in a neighborhood of a point $x_0\in \Omega$, then $\chi(\cdot, \xi)$ 
mimics the characteristic of $u$ of normal direction $\xi=(\xi_1, \xi_2)$ (see Figure \ref{chara}); formally, if $\xi^\perp=(-\xi_2, \xi_1)=\pm u(x_0)$, 
then either $\nabla \chi(\cdot, \xi)$ locally vanishes (if $u$ is constant in a neighborhood of $x_0$), 
or $\nabla \chi(\cdot, \xi)$ is a measure concentrated on the characteristic $\{X(t, x_0)\}_t$ given by \eqref{eq:syst_dyn} with constant 
measure density $\pm \xi$.
In other words, we have the following ``kinetic formulation" of the problem (see e.g., DeSimone-Kohn-M\"uller-Otto\cite{DKMO01} or 
Jabin-Perthame\cite{Jabin-Perthame}): 
\begin{pro}[Kinetic formulation in $2D$]
\label{pro:smooth}
Let $\Omega\subset \RR^2$ be an open set and $u:\Omega\to \RR^2$ be a smooth vector field. If $u$ satisfies \eqref{eq:contraint2} then 
\be
\label{eq:kinet2D}
\xi^\perp \cdot \nabla_x \chi(\cdot, \xi)=0 \quad \textrm{distributionally in } \quad \Omega  \quad \textrm{ for every } \xi \in \Ss^1. \ee
\end{pro}

We mention that the kinetic formulation \eqref{eq:kinet2D} holds under the weaker Sobolev regularity $W^{1/p, p}$ for $p\in [1,3]$ 
(see  Ignat \cite{Ignat_Cras, Ignat_JFA, confluentes}, DeLellis-Ignat \cite{DelIgn}). Note that the knowledge of $\chi(\cdot, \xi)$ in every direction 
$\xi\in \Ss^1$ determines completely a vector field $u$ with $|u|=1$ due to the averaging formula 
\be\label{eq:aver_form}
u(x)=\frac 1 2 \int_{\Ss^1} \xi \chi(x, \xi)\, d\h^1(\xi)  \quad \textrm{ for a.e. }\, \, x\in \Omega.
\ee
Thanks to \eqref{eq:aver_form}, we deduce that the kinetic formulation \eqref{eq:kinet2D} incorporates the fact that $\nabla \times u=0$ 
(see Proposition \ref{pro:curl} below). Therefore, the curl free condition will be no longer mentioned in the
following statements whenever \eqref{eq:kinet2D} is assumed to hold true for unit length vector fields $u$.

The main question is whether the kinetic formulation \eqref{eq:kinet2D} characterizes the vortex vector field in $\R^2$. First of all, 
the equation \eqref{eq:kinet2D} induces a regularizing effect for Lebesgue measurable unit-length vector fields $u$. Indeed,
classical ``kinetic averaging lemma'' (see e.g. Golse-Lions-Perthame-Sentis \cite{Golse}) shows that a measurable vector-field 
$u:\Omega\to \Ss^1$ satisfying \eqref{eq:kinet2D} belongs to $H^{1/2}_{loc}(\Omega)$ due to the averaging formula \eqref{eq:aver_form}.\footnote{For 
the improved regularizing effect for scalar conservation laws, see Otto~\cite{FO} and Golse-Perthame \cite{GolPer}.} 
Moreover, Jabin-Otto-Perthame \cite{JOP02} improved the regularizing effect by showing that $u$ is locally Lipschitz away from 
vortex point-singularities \footnote{This regularity is optimal, see e.g. Proposition 1 in Ignat \cite{Ignat_JFA}.} and $u$ coincides with the 
vortex vector field around these singularities:

\begin{thm} [Jabin-Otto-Perthame \cite{JOP02}]
\label{thm:JOP}
Let $\Omega\subset \RR^2$ be an open set and $u:\Omega\to \RR^2$ be a Lebesgue measurable vector field satisfying 
$|u|=1$ a.e. in $\Omega$ together with the kinetic formulation \eqref{eq:kinet2D}. Then $u$ is locally Lipschitz continuous inside 
$\Omega$ except at a locally finite number of singular points. Moreover, every singular point $P$ of $u$ corresponds to a vortex 
singularity of topological degree one of $u$, i.e., there exists a sign 
$\gamma=\pm 1$ such that $$u(x)=\gamma \us(x-P) \, \textrm{ for every $x\neq P$ in any convex neighborhood of $P$ in $\Omega$.}$$ 
In particular, if $\Omega=\RR^2$ and $u$ is nonconstant, then $u$ coincides with $\us$ or $-\us$ (up to a translation). 
\end{thm}

This result leads to the following interpretation of the kinetic formulation in dimension $N=2$: the equation \eqref{eq:kinet2D} is a 
selection principle for the viscosity solutions of the eikonal equation $|\nabla \psi|=1$ in the sense that the solutions $\psi$ are smooth (more precisely, they belong to the 
Sobolev space $W\loc^{2, \infty}$) away from point-singularities. Clearly, these solutions are induced by the viscosity solutions of the eikonal 
equation under some appropriate boundary condition. Conversely, in the spirit of Caffarelli-Crandall \cite{CafCra}, it is shown by Ignat \cite{Ignat_JFA} 
and De Lellis-Ignat \cite{DelIgn} that for any vector field $u$ satisfying \eqref{eq:contraint2} together with an initial Sobolev regularity
$W^{1/p,p}$, $p\in [1,3]$ (i.e., excluding jump line singularities) then the kinetic formulation \eqref{eq:kinet2D} holds true
and therefore, one obtains the regularizing effect in Theorem \ref{thm:JOP}. 

\begin{rem}
\label{zero-state}
The result of Jabin-Otto-Perthame \cite{JOP02} was motivated by the study of zero-energy states in a line-energy Ginzburg-Landau model in dimension $2$. More precisely, one considers 
the energy functional $E_\eps: H^1(\Omega, \RR^2)\to \RR_+$ defined for $\eps>0$ as
\be
\label{eq:energy}
E_\eps(u_\eps)=\eps \int_\Omega |\nabla u_\eps|^2 \, dx+\frac 1 \eps \int_\Omega (1-|u_\eps|^2)^2\, dx+ \frac 1 \eps \|\nabla\times u_\eps\|_{H^{-1}(\Omega)}^2, \quad u_\eps \in  H^1(\Omega, \RR^2),\ee
where $\Omega$ is a domain in $\R^2$ and $H^{-1}(\Omega)$ is the dual of the Sobolev space $H^1_0(\Omega)$.
(We refer to \cite{AdLM99, AG99, DKMO01, JOP02, Jabin-Perthame, JK00, RS01} for the analysis of this model.) 
A vector field $u:\Omega\to \RR^2$ is called zero-energy state if there exists a family $\{u_\eps \in  H^1(\Omega, \RR^2)\}_{\eps\to 0}$ satisfying
$$u_\eps\to u \quad \textrm{ in } L^1(\Omega) \quad \textrm{and} \quad E_\eps(u_\eps)\to 0 \quad \textrm{ as } \eps\to 0.$$
Obviously, a zero-energy state $u$ satisfies \eqref{eq:contraint2}. The result of Jabin-Otto-Perthame \cite{JOP02} 
shows that every zero-energy state $u$ satisfies \eqref{eq:kinet2D} and therefore, $u$ shares the structure stated in Theorem \ref{thm:JOP}. 
\end{rem}

\subsection{Kinetic formulation in dimension $N\geq 3$}
Our main interest consists in defining a kinetic formulation for the vortex vector field in dimension $N\geq 3$. Let 
$\Omega\subset \RR^N$ be an open set and $u:\Omega\to \RR^N$ be a Lebesgue measurable vector field. For every direction 
$\xi\in \Ss^{N-1}$, we consider the characteristic function $\chi(\cdot, \xi)$ defined at \eqref{eq:weak_cara} and we denote 
the orthogonal hyperplane to $\xi$ by
\[
 \xi^\perp:=\{v\in \RR^N\, :\, v\cdot \xi=0\}.
\]
\begin{df}[Kinetic formulation]
We say that a measurable vector field $u$ satisfies the kinetic formulation if the following equation holds true: 
\be
\label{eq:kineN}
v \cdot \nabla_x \chi(\cdot, \xi)=0 \quad \textrm{distributionally in } \, \Omega \, \textrm{ for every $\xi \in \Ss^{N-1}$ and 
$v\in \xi^\perp$}.\ee
\end{df}

Roughly speaking, \eqref{eq:kineN} means that $\nabla_x \chi(\cdot, \xi)$ is a distribution pointing in direction $\pm \xi$.
Note that the kinetic formulation \eqref{eq:kineN} only carries out the information of the direction of the vector field $u$ (i.e., it gives no information of the
Euclidean norm of $u$). 
Imposing the unit-length constraint, $u$ will satisfy a similar averaging formula \eqref{eq:aver_form} which justifies that the 
curl-free constraint $\nabla \times u=0$ is incorporated in the kinetic formulation \eqref{eq:kineN}.
\begin{pro}
\label{pro:curl}
Let $N\geq 2$, $\Omega\subset \RR^N$ be an open set and $u:\Omega\to \RR^N$ be Lebesgue measurable with $|u|=1$ a.e. in $\Omega$. Then 
\be
\label{eq:averN}
u(x)=\frac1{V_{N-1}} \int_{\Ss^{N-1}} \xi \chi(x, \xi)\, d\h^{N-1}(\xi)  \quad \textrm{ for a.e. }\, \, x\in \Omega,\ee
where $V_{N-1}$ is the volume of the unit ball in $\R^{N-1}$.
Moreover, if $u$ satisfies the kinetic formulation \eqref{eq:kineN} then $\nabla \times u=0$ distributionally in $\Omega$.
\end{pro}

\begin{rem}
\label{rem:contre}
We highlight that Proposition \ref{pro:smooth} is \textit{false} in dimension $N\geq 3$, i.e., there are smooth curl-free vector fields with values into the unit sphere $\Ss^{N-1}$ that do not satisfy the kinetic formulation \eqref{eq:kineN}. 
 For example, in dimension $N=3$, considering the vortex-line vector field 
 $$u_0(x)=\frac{(x_1,x_2,0)}{\sqrt{x_1^2+x_2^2}} \quad \textrm{in} \quad \Omega=\{x=(x_1, x_2, x_3) \in \RR^3 : x_2>1 \},$$ 
 then $u_0$ is smooth in $\Omega$ and satisfies \eqref{eq:contraint2}. However, \eqref{eq:kineN} fails. Indeed, let $\xi=\frac{1}{\sqrt{2}}(1,0,1)$. Then 
 $u_0(x) \cdot \xi=0$ for $x\in \Omega$ is equivalent with $x_1=0$ and therefore,
 \[
 \nabla_x \chi(.,\xi) = e_1 \mathcal{H}^2 \llcorner \{x\in \Omega\, :\, x_1=0\},
 \]
where $e_1=(1,0,0)$. Now, taking $v=\frac{1}{\sqrt{2}}(-1,0,1)$, we have $v \cdot \xi=0$ (i.e., $v \in \xi^{\perp}$) and $v \cdot \nabla_x \chi(.,\xi) \neq 0$ in $\mathcal{D}'(\Omega)$.
\end{rem}

\medskip

As Remark \ref{rem:contre} has already revealed, the kinetic equation \eqref{eq:kineN} in dimension $N\geq 3$ plays a different role than in dimension $N=2$ because the gradient $\nabla \chi(\cdot, \xi)$ is expected to concentrate on hypersurfaces (not on the line characteristics of $u$). In fact, the geometric 
interpretation of \eqref{eq:kineN} can be regarded in terms of the stream function $\psi$ of a nonconstant vector field $u=\nabla \psi$:
the level sets of $\psi$ are expected to be pieces of spheres of codimension one where the characteristics of $u$ represent the normal directions 
to these spheres.

\begin{thm}\label{pro:level}
Let $N\geq 3$, $\Omega\subset \RR^N$ be an open set and $\psi:\Omega\to \RR$ be a smooth stream function such 
that $u=\nabla \psi$ satisfies the kinetic formulation \eqref{eq:kineN}. 
Assume that $|u|$ never vanishes on a level set $\{x\in \Om\,:\, \psi(x)=\alpha\}$ for some $\alpha\in \RR$ 
and let $\Se$ be a connected component of  $\{\psi=\alpha\}$. 
Then $\Se$ is locally a totally umbilical hypersurface, that is either a piece of a $N-1$ sphere or a piece of a hyperplane.
\end{thm}

Note that Theorem \ref{pro:level} fails in dimension $N=2$: a level set of a smooth stream function $\psi$ of $u=\nabla \psi$ satisfying \eqref{eq:contraint2} (and therefore, $u$ satisfies the kinetic formulation \eqref{eq:kinet2D} by Proposition \ref{pro:smooth}) does not have in general constant curvature.\footnote{If $\Gamma\subset \R^2$ is a smooth curve of nonconstant curvature, then one takes $\psi$ to be the distance function to $\Gamma$ in a small neighborhood $\Omega$ of $\Gamma$ (with the convention that $\Gamma$ is withdrawn from that neighborhood, i.e., $\Gamma\cap \Omega=\emptyset$, so that $\psi$ is smooth in $\Omega$). } 

\section{Main results}

Our main result shows that the kinetic formulation \eqref{eq:kineN} is a characterization of the vortex vector field $\us$ in dimension 
$N \geq 3$.

\begin{thm}
\label{thm:main}
Let $N\geq 3$, $\Omega\subset \RR^N$ be a connected open set and $u:\Omega\to \RR^N$ be a nonconstant Lebesgue
measurable vector field satisfying $|u|=1$ a.e. in $\Omega$ together with the kinetic equation \eqref{eq:kineN}. 
Then $u$ coincides with the vortex vector field $\us$ or $-\us$ up to a translation.
\end{thm}

Note that in dimension $N=2$, this result is true for the domain $\Omega=\RR^2$, but it is in general false for other domains $\Omega$ 
where there exist nonconstant smooth vector fields $u$ in $\Omega$ different than vortex vector fields that satisfy \eqref{eq:contraint2} 
and thus, \eqref{eq:kinet2D} (by Proposition \ref{pro:smooth}). The main difference in dimension $N\geq 3$ is the following: 
if $u$ is a smooth vector field with \eqref{eq:contraint2} that is neither constant nor a vortex vector field, then the kinetic 
formulation \eqref{eq:kineN} doesn't hold for $u$ (see Remark \ref{rem:contre}). Hence, in dimension $N\geq 3$, the zero-energy states of $E_\eps$ defined at \eqref{eq:energy}
does not satisfy in general the kinetic equation \eqref{eq:kineN}. Therefore, the kinetic formulation \eqref{eq:kineN} is more rigid in 
dimension $N\geq 3$ since it selects only the vortex vector fields as they correspond to smooth solutions of the eikonal 
equation with level sets of constant sectional curvature (by Theorem \ref{pro:level}).

Let us explain the strategy of the proof of Theorem \ref{thm:main}. The key point lies on a relation of order of the level sets of
the stream function associated to $u$: for every two Lebesgue points $x, y\in \Omega$ of $u$ such that the segment 
$[x,y]\subset \Omega$ and for every direction $\xi\in \Ss^{N-1}$ orthogonal to $x-y$, one has
$$u(x)\cdot \xi>0 \Rightarrow u(y)\cdot \xi\geq 0.$$ The next step consists in defining the trace of $u$ on each segment 
$\Sigma\subset \Omega$; more precisely, 
similar to the procedure of Jabin-Otto-Perthame \cite{JOP02}, there exists a trace
$\tilde u\in L^\infty(\Sigma, \Ss^{N-1})$ of $u$ such that  $u(P)=\tilde u(P)$ for each Lebesgue point $P \in \Sigma$ of $u$. Moreover, if the trace
$\tilde u$ of $u$ is collinear with the segment $\Sigma$ at some Lebesgue point, then $\tilde u$ is $\mathcal{H}^1$-almost everywhere collinear with 
$\Sigma$ (which coincides with the classical principle of characteristics for smooth vector fields $u$). The final step consists 
in proving that every two characteristics are coplanar. Then the conclusion follows by the following geometrical fact specific to dimension $N\geq 3$:

\begin{pro}\label{prop:geometry}
Let $N\geq 3$ and ${\cal D}$ be a set of lines in $\RR^N$ such that every two lines 
of $\cal D$ are coplanar, but ${\cal D}$ is not planar (i.e., there is no two-dimesional plane containing $\cal D$). Then either all lines of $\cal D$ are collinear, or all lines of $\cal D$ pass through a 
same point (that is a vortex point). 
\end{pro}

In view of Theorem \ref{thm:main}, it is natural to ask if one can characterize other type of unit-length curl-free vector fields $u$
by weakening the kinetic formulation \eqref{eq:kineN}, in particular, vector fields having a vortex-line singularity. In dimension $N\geq 3$, the prototype of a vortex-line vector field is given by 
$$u_0(x',x_N)=\nabla |x'|$$ where $x=(x', x_N)$, $x'=(x_1,\ldots,x_{N-1})$; clearly, $u_0$ is smooth away from the vortex-line 
$\{x\in \R^N\,:\, x'=0\}$ where \eqref{eq:contraint2} holds true. Denoting $$\mathcal{E}:=\{\xi \in \Ss^{N-1} \,: \, \xi_N=0 \}=
\Ss^{N-2}\times \{0\},$$ 
within the notation \eqref{eq:weak_cara}, we have that $u_0$ satisfies the following kinetic formulation in $\Omega=\R^N$:
\begin{equation}\label{eq:kineNW}
\forall \xi \in \mathcal{E}, \quad \forall v \in \xi^{\perp}, \quad v \cdot \nabla_x \chi(.,\xi)=0 \quad \text{ in } \quad \mathcal{D}'(\Omega).
\end{equation}
Note that \eqref{eq:kineNW} is a weakened form of \eqref{eq:kineN}: the quantity $v \cdot \nabla_x \chi(.,\xi)$ vanishes
 for directions $\xi \in \mathcal{E}$ (and $v\in \xi^{\perp}$) and fails to vanish for $\h^{N-1}$-a.e. direction $\xi \in \Ss^{N-1}$. As opposed to \eqref{eq:kineN} (in view of \eqref{eq:averN}), the kinetic formulation \eqref{eq:kineNW} 
does not force a unit-length vector field $u$ to be curl-free; it only implies that $$\nabla'\times \frac{u'}{|u'|}=0  \quad \text{ in } \quad \{|u'|\neq 0\}=\{u\neq \pm e_N\}$$
where $e_N=(0, \dots, 0, 1)$, $u'=(u_1, \dots, u_{N-1})$ and $\nabla'=(\partial_1, \dots, \partial_{N-1})$. Since we look for a characterization of vortex-line vector fields (that are in particular curl-free), we will impose that
\be
\label{differ_deriv}
\partial_k u_N=\partial_N u_k \quad \text{ in } \quad \Omega,  \text{ for } \, k=1, \dots, N-1. 
 \ee 

We will prove the following result:

\begin{thm}
\label{thm:vortex-line}
Let $N\geq 4$, $\Omega\subset \R^N$ be an open set and $u : \Omega \to \R^N$ be a Lebesgue measurable vector field satisfying
$|u|=1$ a.e. on $\Omega$ together with \eqref{eq:kineNW} and \eqref{differ_deriv}. Then in every ball included in $\{x\in \Omega\, :\, u(x)\neq \pm e_N\}$, there exists a stream function $\psi=\psi(\alpha, \beta)$ solving the eikonal equation in dimension $2$ such that 
$$u(x)=\nabla_x[ \psi(\alpha, \beta)]$$
where
\begin{enumerate}
\item[1)] either $\alpha=|x'-P'|$ and $\beta=x_N$ for some point $P'\in \R^{N-1}$;

\item[2)] or $\alpha=w'\cdot x'$ and $\beta=x_N$ for some vector $w'\in \Ss^{N-2}$.
\end{enumerate}
\end{thm}

Therefore, the weakened kinetic formulation \eqref{eq:kineNW} (together with \eqref{differ_deriv}) is not enough to select vortex-line vector fields which correspond to the stream function $\psi(\alpha, \beta)=\pm \alpha$ in the case $1)$ of Theorem \ref{thm:vortex-line}. Similar results to Theorem \ref{thm:vortex-line} hold for similar kinetic formulations corresponding to vector fields having vortex-sheets singularities of dimension $k$ in $\R^{N}$ with $N\geq k+3$.

\bigskip

The outline of this paper is as follows: in Section 3, we characterize the level sets of smooth stream functions associated to vector fields that
satisfy the kinetic formulation \eqref{eq:kineN}. In particular, we prove Proposition \ref{pro:smooth} and Theorem \ref{pro:level}. Section 4 is devoted 
to prove fine properties of Lebesgue points of $u$ needed in Section 5 where the notion of trace on lines for a vector field
$u$ satisfying \eqref{eq:kineN} is defined. Section 6 is the core of this paper: using this notion of trace and the geometric arguments of
Proposition \ref{prop:geometry}, we prove our main result in Theorem \ref{thm:main}. 
The last section deals with the
study of the weakened kinetic formulation \eqref{eq:kineNW}.

\section{Level sets of the stream function}
This section is devoted to the study of the level sets of smooth stream functions $\psi$ associated to vector fields $u=\nabla \psi$
satisfying \eqref{eq:kineN}. We start by proving that $|\nabla \psi|$ is locally constant on each level set of $\psi$.

\begin{lem}
\label{lem:derivative}
Let $N\geq 2$, $\Omega\subset \RR^N$ be an open set and $\psi:\Omega\to \RR$ be a smooth stream function such that 
$u=\nabla \psi$ satisfies the kinetic formulation \eqref{eq:kineN}. Assume that $|u|$ never vanishes on a level 
set $\{x\in \Omega\,:\, \psi(x)=\alpha\}$ for some $\alpha\in \RR$ and  let $\Se$ be a connected component of  $\{\psi=\alpha\}$. 
Then $|u|$ is constant on $\Se$. Moreover, there exists a neighborhood $\om$ of $\Se$, a smooth solution $\tilde \psi:\om\to \RR$ 
of the eikonal equation and a diffeomorphism $t \mapsto F(t)$ such that $\psi=F(\tilde \psi)$ in $\om$ (in particular, $\nabla \tilde \psi$ satisfies \eqref{eq:kineN}). 
\end{lem}
\proof{}
Since $|u| \neq 0$ on $\Se$ and $u$ is smooth in $\Omega$, we can define $$v=\frac{u}{|u|} \quad \textrm{ in a neighborhood of $\Se$}.$$ For simplicity of the writing, we suppose that $\Omega$ is this neighborhood, i.e., $|u|\neq 0$ in $\Omega$. Then
$v$ satisfies \eqref{eq:kineN} because $u$ satisfies it, too; $v$ being smooth in $\Omega$, then\footnote{The proof of Proposition \ref{pro:curl} is independent of Lemma \ref{lem:derivative}; we will
admit it here and prove it later in Section 4.} Proposition \ref{pro:curl} implies $\nabla \times v=0$  in $\Omega$. As a consequence, in any simply connected domain $\om \subset \Omega$, the Poincar\'e lemma yields
the existence of a smooth function $\tilde{\psi}$ such that $v=\frac{u}{|u|}=\nabla \tilde{\psi}$ in $\om$, i.e., 
\[
 \nabla \psi =u= |u|v=|u| \nabla \tilde{\psi} \text{ \quad in \quad $\om$}.
\]
Therefore, $\psi$ and $\tilde{\psi}$ have the same level sets in $\om$. W.l.o.g, we may assume that $\tilde{\psi}=0$ on $\om \cap \Se$. Now, for every $P' \in \om \cap \Se$, we consider the flow associated to $v$:
\begin{align}\label{EDO}
\left\{
\begin{array}{l}
 \dot{X}(P',t)=\nabla \tilde{\psi}(X(P',t))\\
 X(P',0)=P'.
\end{array}
\right.
\end{align}
Call $I_{P'}$ the maximal interval where the solution $X(P',.)$ exists. Obviously, the flow is unique and smooth satisfying the following:
\[
 \ddot{X}(P',t)=\nabla^2 \tilde{\psi}(X) \cdot \dot{X}=\nabla^2 \tilde{\psi}(X) \cdot \nabla \tilde{\psi}(X)=0 \text{ \quad in \quad $I_{P'}$}
\]
because $\nabla^2 \tilde{\psi}$ is a symmetric matrix and $|\nabla \tilde{\psi}|=1$ in $\om$. Consequently, $\dot{X}(P',\cdot)$ is constant in $I_{P'}$ so that $\nabla \tilde{\psi}(X(P',t))=\nabla \tilde{\psi}(P')$, $\frac{d}{dt}[\tilde \psi(X(P',t))]=1$ and $X(P',t)=P'+t \nabla \tilde{\psi}(P')$. Therefore, since $\tilde{\psi}=0$ on $\om \cap \Se$, we have:
\[
 \tilde{\psi}(X(P',t))=t \quad \textrm{for all $P' \in \om \cap \Se$ and 
$t \in I_{P'}$.}
\]
Identifying the level sets of $\tilde{\psi}$ (and of $\psi$, too) using the flow, i.e., $\{\tilde{\psi}=t \}=\{X(P',t)\,:\,  P' \in \om \cap \Se \}$, we can define $$F(t):=\psi(X(P',t)), \quad \textrm{for}\quad  P' \in \om \cap \Se, t\in I_{P'}.$$ 
The function $F$ is a diffeomorphism: $F$ is smooth (because $\psi$ and $X$ are smooth, too) and we have 
\[
 \frac{d}{dt}F(t)=\nabla \psi(X(P',t)) \cdot \dot{X}(P',t)\stackrel{\eqref{EDO}}{=} \nabla \psi(X(P',t)) \cdot \frac{\nabla \psi}{|\nabla \psi|}(X(P',t))=|u|(X(P',t))\neq 0.
\]
In particular, $|u|$ is constant on $\{ \tilde{\psi}=0 \}=\{\psi=F(0) \}=\om\cap \Se$. Since $\om$ was arbitrarily chosen, we deduce that $|u|$ is locally constant
on $\Se$; because $\Se$ is connected, it follows that $|u|$ is constant on $\Se$. Since the flow $\{X(P',t)\, :\, P'\in \Se, t\in I_{P'}\}$ covers a neighborhood of $\Se$, the last statement of the lemma follows, too.  \qed

\subsection{The case of dimension $N=2$}
In the special case of dimension $N=2$, we start by proving that every smooth curl-free vector field of unit length
satisfies the kinetic formulation \eqref{eq:kinet2D}. This result can be found already in the works of DeSimone-Kohn-M\"uller-Otto\cite{DKMO01} or 
Jabin-Perthame\cite{Jabin-Perthame}. For completeness of the paper, we will present two easy and self-contained proofs. The first one is based on the geometry of the flow \eqref{eq:syst_dyn} (as heuristically exposed at Section \ref{chap_div}), while the second proof is based on the concept of entropy introduced in \cite{DKMO01}.

\proof{ of Proposition \ref{pro:smooth}: First method}
We can assume that $\xi=e_1$ and $\xi^{\perp}=e_2$ (otherwise, one considers a rotation $R\in SO(2)$ such that $e_1=R \xi$ and 
$\tilde{u}(x):=Ru(R^{-1}x)$ in a neighborhood of a point $x\in \Omega$).
Naturally, $\Omega$ can be written as a countable reunion of squares whose edges are parallel with $e_1$ and $e_2$. Therefore, using a partition of unity,
it is enough to prove the statement for $\Omega=(-1,1)^2$:
\[
\forall \varphi \in C_c^{\infty}(\Omega), \quad 0= \int_{\Omega} \varphi \, \xi^\perp\cdot \nabla_x \chi(x,\xi) \diff x\stackrel{\xi=e_1}{=}\int_{\Omega} \varphi \partial_2 \chi(x,e_1) \diff x = -\int_{\Omega \cap \{u_1>0\}} \partial_2 \varphi \diff x.
\]
For that, we consider the flow \eqref{eq:syst_dyn} and by the proof of Lemma \ref{lem:derivative}, we have that for every $x\in \Om$, 
$\{X(t,x) \}_t$ is a straight line given by $X(t,x)=x+tu(x)$ and $u(X(t,x))=u(x)$ for all $t$. 
Since $u$ is smooth, there is no crossing  between two characteristics in $\Omega$. We claim that:
\[
 \Omega \cap \{u_1>0\}= \bigsqcup_{k\in K} A_k,
\]
where $\{A_k\}_{k \in K}$ is a (at most) countable disjoint set of rectangles of type $(a_k,b_k) \times (-1,1)\subset \Om=(-1,1)^2$. Indeed, if 
$x \in \Omega \cap \partial \{u_1>0 \}$ then $u_1(x)=0$ and $u(x) \parallel e_2$: therefore, for all $t$, $X(t,x) \parallel e_2$. So 
the set $\Omega \cap \partial \{u_1>0\}$ is a (at most) countable set of vertical segments $\{x_1\}\times (-1,1)$ inside $\Omega$ with $x_1\in \{a_k, b_k\}_{k\in K}\subset [-1,1]$, and the claim is proved. Now,
for $\varphi \in C_c^{\infty}(\Omega)$,
\[
 \int_{\Omega \cap \{u_1>0\}} \partial_2 \varphi =\sum_k \int_{A_k} \partial_2 \varphi=\sum_k \int_{a_k}^{b_k} \int_{-1}^1 \partial_2 \varphi= 0,
\]
because $\partial_2 \varphi$ can be seen as a signed Radon measure for $\varphi \in C_c^{\infty}(\Omega)$ and the proposition is proved. \qed

\proof{ of Proposition \ref{pro:smooth}: Second method}
The following proof links the kinetic formulation
\eqref{eq:kinet2D} with the theory of entropy solutions for scalar conservation laws (see e.g., \cite{DKMO01}). Indeed, if $u$ is a smooth vector field
satisfying \eqref{eq:contraint2}, then formally, $u_1=-h(u_2):=\pm \sqrt{1-u_2^2}$ so that $\nabla \times u=0$ can be rewritten as:
\begin{align}\label{conservation_law}
 \partial_1 u_2 + \partial_2 [h(u_2)]=0;
\end{align}
thus, $u_2$ can be formally interpreted as a solution of the above scalar conservation law in the variables $(time, space)=(x_1, x_2)$.
Based on the concept of entropy solution of \eqref{conservation_law} introduced via the pairs (entropy, entropy-flux), the following applications (called ``elementary entropies") were used in \cite{DKMO01}. More precisely, for every $\xi\in \Ss^1$, $\Phi^\xi:\Ss^1\to \RR^2$ is defined as 
$$\textrm{for $z \in \Ss^1$,} \quad \Phi^\xi(z)=\begin{cases} \xi^\perp &\quad \textrm{ for }\, \, z\cdot \xi>0, \\
0 &\quad \textrm{ for }\, \, z\cdot \xi\leq 0.
\end{cases}$$ 
Then the kinetic formulation \eqref{eq:kinet2D} writes as
\be
\label{eq:zero_entropy_prod}
\nabla\cdot \big[\Phi^\xi(u)\big]=0 \quad \textrm{distributionally in} \; \Omega.
\ee
In order to prove \eqref{eq:zero_entropy_prod}, we will approximate $\Phi^\xi$ by a sequence of smooth maps $\{\Phi_k:\Ss^1\to \RR^2\}$ such that
$\{\Phi_k\}$ is uniformly bounded, $\lim_k \Phi_k(z)=\Phi^\xi(z)$ for every $z\in \Ss^1$ and $\Phi_k$ satisfies \eqref{eq:zero_entropy_prod} for every $k$. 
Following the ideas in \cite{DKMO01} (see also \cite{IMpre}), this smoothing result comes from the following observation: there exists a (unique) 
$2\pi$-periodic piecewise $C^1$ function $\f: \RR \to \RR$ associated to $\Phi^\xi$ via the equation
\be
\label{eq:entro_decomp}
\Phi^\xi(z)=-\f'(\theta)z+\f(\theta)z^\perp \quad \textrm{for every } z=e^{i\theta}\in \Ss^1.
\ee
In fact, $\f$ is given by:
$$\f(\theta)=\Phi^\xi(z)\cdot z^\perp=\xi\cdot z {\mathds{1}}_{\{z\cdot\xi>0\}}=\cos(\theta-\theta_0) 
{\mathds{1}}_{\{\theta-\theta_0\in (-\pi/2, \pi/2)\}} \quad \textrm{ for }\, \, z=e^{i\theta}, \, \theta\in (-\pi+\theta_0, \pi+\theta_0),$$
where $\xi=e^{i\theta_0}\in \Ss^1$ with $\theta_0\in(-\pi, \pi]$.  
In \eqref{eq:entro_decomp}, the distributional derivative $\f'$ is given by 
$$\f'(\theta)=-\sin (\theta-\theta_0) {\mathds{1}}_{\{\theta-\theta_0\in (-\pi/2, \pi/2)\}} \quad \textrm{ for}\, \, 
\theta\in (-\pi+\theta_0, \pi+\theta_0).$$
Now, one regularizes $\f$ by $2\pi-$periodic functions $\f_k\in C^\infty(\RR)$ that are uniformly bounded in
$W^{1, \infty}(\RR)$ and $\lim_k \f_k(\theta)=\f(\theta)$ as well as
 $\lim_k \f'_k(\theta)=\f'(\theta)$ for every $\theta\in \RR$. Then we define $\Phi_k$ as in \eqref{eq:entro_decomp} for the functions $\f_k$:
 $$\Phi_k(z)=-\f'_k(\theta)z+\f_k(\theta)z^{\perp} \, \quad \textrm{for $z=e^{i \theta}\in \Ss^1$}.$$ 
Let us now check that  $\{\Phi_k\}_k$ are indeed the desired (smooth)
 approximating maps of $\Phi^\xi$. For that, first, note that differentiating the above equation defining $\Phi_k$, one obtains that 
 \be
 \label{eq:deriv_angul}
 \frac{\partial \Phi_k}{\partial \theta}(z)\cdot z^\perp=0 \quad \textrm{ for every }\quad z=e^{i\theta}\in \Ss^1.\ee
 Next, we prove that $\Phi_k$ satisfies \eqref{eq:zero_entropy_prod}. Indeed, we can locally write $u=e^{i\Theta}$ in every ball $B\subset \Omega$ 
 for some smooth lifting $\Theta:B\to \RR$ that satisfies
 $$\nabla \Theta\cdot u=\nabla\times u=0 \quad \textrm{in}\quad B.$$ This means that $\nabla \Theta=\lambda u^\perp$ in 
 $B$ for some smooth function $\lambda:B\to \RR$. Therefore, it follows
$$\nabla \cdot \big[\Phi_k(u)\big]=\frac{\partial \Phi_k}{\partial \theta}(e^{i\Theta})
\cdot \nabla \Theta=\lambda \frac{\partial \Phi_k}{\partial \theta}(u)\cdot 
u^\perp\stackrel{\eqref{eq:deriv_angul}}{=}0 \quad \textrm{in}\quad B.$$ 
Passing at the limit $k\to \infty$, the dominated convergence theorem yields: 
$$\int_B \Phi^\xi(u)\cdot \nabla \zeta\, dx=0 \quad
\textrm{for every $\zeta\in C^\infty_c(B)$}.$$
The conclusion is now straightforward.
\qed

\medskip

Note that another interest of this second method is that it can be adapted to vector field $u \in W^{\frac{1}{p},p}$ for $p\in [1,3]$. 
For such vector fields, there is a priori no trace of $u$ on a segment so that the flow \eqref{eq:syst_dyn} does not have a proper meaning anymore;
see \cite{Ignat_JFA} and \cite{DelIgn} for more details.

\subsection{The case of dimension $N\geq 3$}

The aim of this subsection is to prove Theorem \ref{pro:level}. We divide the proof in several steps, each step being stated as a lemma.

\begin{lem}\label{lem:represent_surf}
Let $\Omega\subset \R^N$ be an open set and $u : \Omega \to \R^N$ be a smooth vector field satisfying \eqref{eq:kineN}. We denote by
$$\tilde{\Omega}:= \{x \in \Omega \, : \, u(x)\neq 0,\,  \nabla \big(\frac{u}{|u|}\big)(x)\neq 0\}$$
and for every $x\in \tilde{\Omega}$, 
$$\Ss_x:= u(x)^\perp\cap \Ss^{N-1}=\{\xi \in \Ss^{N-1} : u(x) \cdot \xi=0\} \approx \Ss^{N-2}.$$
Then we have for all $x \in \tilde{\Omega}$ and for $\mathcal{H}^{N-2}$-a.e. $\xi \in \Ss_x$ that the set
\[
 \{y \in \tilde \Omega : u(y)\cdot \xi=0\}=\tilde \Omega\cap \partial \{u \cdot \xi>0\}
\]
is a hyperplane around $x$ that is oriented by the normal vector $\xi$. Moreover,
\be\label{chi_representation}
\nabla_x \chi(.,\xi)=\pm \xi \mathcal{H}^{N-1} \llcorner\partial\{u \cdot \xi>0\} \quad \text{ locally around $x$}.
\ee
\end{lem}
\proof{}
As in the proof of Lemma \ref{lem:derivative}, we set $v=\frac{u}{|u|}$ on $\tilde \Om$. Then $v$ is a smooth unit-length vector field in $\tilde \Omega$ that satisfies \eqref{eq:kineN} (because $u$ satisfies it, too) and by Proposition \ref{pro:curl}, we have that $v$ is curl-free in $\tilde \Omega$.
Let $x \in \tilde{\Omega}$, in particular, $\nabla v(x)\neq 0$. First, we show that $\{y \in \tilde \Omega : u(y) \cdot \xi =0 \}$ is a smooth $N-1$ manifold around $x$. Since $v$ is curl-free, we know that $\nabla v(x)=(\partial_j v_i(x))_{i,j}$ is symmetric. By differentiating the relation $|v(x)|=1$, it follows:
\begin{align*}
\nabla v(x)^T v(x)= \nabla v(x) v(x)=0.
\end{align*}
That means $v(x) \in \text{Ker }\nabla v(x)$.
We will prove that
$$\mathcal{H}^{N-2}(\Ss_x \cap \text{Ker } \nabla v(x))=0.$$
Assume by contradiction that $\Ss_x \cap \text{Ker } \nabla v(x)$ has positive $\h^{N-2}$-measure. Since $\text{Ker } \nabla v(x)$ is a linear space, then one would have that $\Ss_x \subset \text{Ker }\nabla v(x)$, i.e., $\nabla v(x) \xi =0$ for all $\xi \in \Ss_x$. Moreover, since $v(x) \in \text{Ker }\nabla v(x)$ and $\Ss_x\subset v(x)^\perp$, it follows that $\nabla v(x)=0$ which is a contradiction with the assumption $\nabla v(x)\neq 0$.
Therefore, $\nabla v(x) \xi \neq 0$ for $\mathcal{H}^{N-2}$-a.e. $\xi \in \Ss_x$ and $\{y \in \tilde \Omega : v(y) \cdot \xi =0 \}=\{y \in \tilde \Omega : 
u(y) \cdot \xi =0 \}$ is a smooth $N-1$ manifold around $x$.


 
It remains to prove that this manifold is a piece of hyperplane oriented by $\xi$ where
\eqref{chi_representation} holds true. For that, set $\varphi \in C_c^{\infty}(\tilde \Omega,\R^N)$ be supported in a ball $B\subset \tilde \Omega$ centered at $x$. 
By the Gauss theorem, we have
\begin{align*}
-\langle \nabla_x \chi(.,\xi),\varphi \rangle=\int_{B} \nabla \cdot \varphi(y) \chi(y,\xi) \diff y  &= \int_{  \{y\in B\, :\, u(y)\cdot \xi>0\} } \nabla \cdot \varphi \diff y\\
&=\int_{B\cap \partial \{u\cdot \xi>0\}} \varphi \cdot \nu \diff \mathcal{H}^{N-1}(y)
\end{align*}
where $\nu$ is the unit outer normal vector at the $N-1$ manifold $\partial \{ u(y) \cdot \xi>0 \}$. This proves that locally around $x$, we have
$$\nabla_x \chi(x,\xi)=-\nu \mathcal{H}^{N-1} \llcorner \bigg(B\cap \partial \{u\cdot \xi >0\}\bigg).$$ Because of \eqref{eq:kineN}, we know that
$\nabla_x \chi(x,\xi)$ and $\xi$ are collinear. Since $\nu$ is smooth on $B\cap \partial \{u \cdot \xi>0 \}$, this implies
$\nu=\xi$ or $\nu=-\xi$ on $B\cap \partial \{u\cdot \xi >0\}$. The conclusion is now straightforward. \qed 
\bigskip

We now state the following result which is the key point in proving Theorem \ref{pro:level}.
\begin{lem}\label{lem:second_form_propid}
Under the hypothesis of Theorem \ref{pro:level}, every point $x \in \Se$ is an umbilical point, i.e., 
there exists $\lambda(x)\in \R$ such that:
\[
Du(x)=\lambda(x) Id \colon T_x \Se  \longrightarrow  \RR^{N-1} 
\]
where $u$ is proportional with the Gauss map on $\Se$, $T_x \Se$ is the tangent plane at the hypersurface $\Se$ at $x$ and $Id$ is the identity matrix. 
\end{lem}
\proof{} Recall that $|u|$ is constant on $\Se$ by Lemma \ref{lem:derivative} so that $u/|u|$ is the normal vector (i.e., the Gauss map) at the hypersurface $\Se$. Therefore,
$$D\big(\frac{u}{|u|}\big|_{\Se}\big)=\frac{1}{|u|}D\big(u\big|_{\Se}\big) \quad \textrm{ in }\quad \Se,$$
where $D(u\big|_{\Se})$ is the differential of $u$ restricted to $\Se$ as a map with values into the sphere $\Ss^{N-1}$ (up to the multiplicative constant $|u|$). 
As in the proof of Lemmas \ref{lem:derivative} and \ref{lem:represent_surf}, we may assume that $u$ never vanishes in $\Omega$ and set
$v=\frac{u}{|u|}$ in $\Omega$. Then $v$ is a smooth unit-length vector field in $\Omega$ that satisfies \eqref{eq:kineN} 
and by Proposition~\ref{pro:curl}, $v$ is curl-free so that locally $v=\nabla \tilde \psi$ for a smooth stream function $\tilde \psi$. Since $\nabla \psi=u=|u|\nabla \tilde \psi$, we know that $\psi$ and $\tilde \psi$ have the same level sets, in particular, $\Se$ is a level set of $\tilde \psi$. Therefore, replacing $u$ by $v$, we may assume in the following that 
$$|u|=1\quad \textrm{ in }\quad \Omega.$$
Let $x\in \Se$. We want to show that $x$ is an umbilical point of $\Se$. This is clear if $\nabla u(x)=0$. Therefore, we assume in the following that $x\in \tilde \Omega\cap \Se$ defined at Lemma \ref{lem:represent_surf}, i.e., $$\nabla u(x)\neq 0.$$ 
 Since \eqref{eq:averN} holds for the unit-length vector field $u$, we obtain by differentiating \eqref{eq:averN}:
\[
\nabla u(x)=\frac1{V_{N-1}} \int_{\Ss^{N-1}} \xi \otimes \nabla_x \chi(x,\xi), 
\]
where $V_{N-1}$ is the volume of the unit ball in $\R^{N-1}$. The above integrant is to be understood as an absolutely continuous measure with respect to the Hausdorff $\h^{N-2}$ measure concentrated on the set $\Ss_x$ (defined at Lemma \ref{lem:represent_surf}). For that, we check first that the support of the integrand lies on $\Ss_x$. Indeed, if $\xi \in \Ss^{N-1}$ with $u(x) \cdot \xi \neq 0$, then $\nabla_x \chi(\cdot,\xi)=0$ in the open set $\{ u \cdot \xi \neq 0\}$ around $x$. Therefore, the integrand has support on the set $\xi \in \Ss_x$ where \eqref{chi_representation} holds true
for $\h^{N-2}$-a.e. $\xi \in \Ss_x$, the density of the measure being equal with $\pm \xi \otimes \xi  \h^{N-2}\llcorner \Ss_x$. Since $\Ss_x\subset u(x)^\perp=T_x \Se$, the density $\xi \otimes \xi$ with $\xi \in \Ss_x$ already identifies $\nabla u(x)\equiv D u(x)$. Next we compute this quantity by exploring the sign of the density $\pm \xi \otimes \xi$:

\medskip

\noindent \textbf{Case $N=3$}. We show that there are at most two nonzero vectors $\pm \xi_0\in \Ss_x\approx \Ss^1$ such that $\nabla u(x) \xi_0=0$.  
Assume by contradiction that there are more than two vectors as above, i.e., there exists another nonzero vector $\tilde{\xi}_0 \neq \pm \xi_0$ in $\Ss_x$ such that 
$\nabla u(x) \xi_0=\nabla u(x) \tilde{\xi}_0=0$. 
Because of $|u|=1$, we know that $\nabla u(x) u(x)=0$. Since the set $\{u(x), \xi_0, \tilde{\xi}_0\}$ spans $\R^3$, it implies $\nabla u(x)=0$ which contradicts the hypothesis $x\in \tilde \Omega$. Therefore, $\nabla u(x) \xi\neq 0$ for every $\xi \in \Ss_x\setminus \{\pm \xi_0\}$ (or for every $\xi \in \Ss_x$ if $\xi_0$ does not exist) and by Lemma~\ref{lem:represent_surf}, $\partial \{ u(y) \cdot \xi >0 \}$ is a smooth surface around $x$ oriented by $\xi$. Let $\mathcal{C}_1$ and 
$\mathcal{C}_2$ be the two connected components of $\Ss_x\setminus \{\pm \xi_0\}$ (convention: $\mathcal{C}_1=\mathcal{C}_2=\Ss_x$ in the case $\nabla u(x) \xi\neq 0$ for every $\xi \in \Ss_x$). For $j=1,2$,
we associate to a point $\xi\in\mathcal{C}_j$ the unit outer normal vector field $\nu(\xi)\in \{\pm \xi\}$ at the plane $\partial \{ u \cdot \xi >0 \}$ around $x$.
Since the map $\xi\in\mathcal{C}_j\to \nu(\xi)$ is smooth (by the implicit function theorem) and $\mathcal{C}_j$ is connected, we deduce that $\nu$ is constant on $\mathcal{C}_j$. Thus it follows that
\[
\pi \nabla u(x)= \gamma_1 \int_{\mathcal{C}_1} \xi \otimes \xi \diff \xi +\gamma_2 \int_{\mathcal{C}_2} \xi \otimes \xi \diff \xi,
\]
with $V_2=\pi$ and $\gamma_{1,2} \in \{\pm 1 \}$ (with the convention that $\gamma_1=\gamma_2=\pm \frac 12$ if $\mathcal{C}_1=\mathcal{C}_2=\Ss_x$).
It remains to show that $\int_{\mathcal{C}_j} \xi \otimes \xi \diff \xi$ is proportional with the identity matrix $Id$, $j=1,2$. Up to a rotation, we can 
suppose that $u(x)=e_3$ and $\mathcal{C}_1=\{\xi \in \Ss^1\times\{0\} : \xi_2>0 \}\approx \{(\cos \theta, \sin \theta)\, :\, \theta\in (0, \pi)\}$. We have
\begin{align*}
\int_{\mathcal{C}_1} \xi \otimes \xi \diff \xi & \approx \int_{0}^{\pi} 
\begin{pmatrix}
\cos^2 \theta & \cos \theta \sin \theta \\
\cos \theta \sin \theta & \sin^2 \theta 
\end{pmatrix}
\diff \theta  = \frac{\pi}{2}\, Id
\end{align*}
(the conclusion follows similarly if $\mathcal{C}_1=\mathcal{C}_2=\Ss_x$).


\medskip

\noindent \textbf{Case $N>3$}. Let $\mathcal{C}=\text{Ker } \nabla u(x) \cap \Ss_x$. We know that $u(x) \in \text{Ker } \nabla u(x) $ and $u(x)$ is orthogonal at $\Ss_x$ which is isomorphic to $\Ss^{N-2}$. Since $\nabla u(x)\neq 0$ (i.e., the dimension $\text{Ker } \nabla u(x)$ is at most $N-1$),  we have two situations (as in the case $N=3$):

$\bullet$ either $\dim \text{Ker } \nabla u(x)=N-1$ leading to $\mathcal{C}$ isomorphic to $\Ss^{N-3}$. In this situation, $\Ss_x\setminus \mathcal{C}$ is the partition of two connected sets $\mathcal{C}_1$ and $\mathcal{C}_2$ that are isomorphic to the half sphere
$$\Ss_+^{N-2}=\{\xi=(\xi_1, \dots, \xi_{N-1}) \in \Ss^{N-2}\, :\,  \xi_1>0 \}.$$ 
The same argument as in the case $N=3$ shows that the sign of the unit outer normal field $\nu(\xi)\in \{\pm \xi\}$ at the hyperplane $\partial \{ u \cdot \xi >0 \}$ is constant when $\xi$ covers $\mathcal{C}_j$, $j=1,2$, so that
\[
V_{N-1}\nabla u(x)= \gamma_1 \int_{\mathcal{C}_1} \xi \otimes \xi \diff \xi +\gamma_2  \int_{\mathcal{C}_2} \xi \otimes \xi \diff \xi,
\]
with $\gamma_1, \gamma_2 \in \{\pm 1\}$.

$\bullet$ or $\dim \text{Ker } \nabla u(x)\leq N-2$ leading to the manifold $\mathcal{C}$ of dimension $\leq N-4$. In other words, $\Ss_x\setminus \mathcal{C}$ is connected and covers a.e. point of $\Ss_x$. The above formula holds for $\mathcal{C}_1=\mathcal{C}_2=\Ss_x$ and 
$\gamma_1=\gamma_2=\pm \frac 1 2$. 

\medskip

We now compute $\nabla u(x)$. For that, we may assume (up to a rotation) that $u(x)=e_N$ and $\mathcal{C}_1=\Ss_+^{N-2}$. 
Since $\Ss_+^{N-2}$ is invariant under the change of coordinate $\xi_d \mapsto -\xi_d$ for some $2\leq d\leq N-1$, we have for every $1\leq j\leq N-1$ with $j\neq d$:
\[
 \int_{\Ss_+^{N-2}} \xi_j \xi_d \diff \xi=-\int_{\Ss_+^{N-2}} \xi_j \xi_d \diff \xi=0,
\]
leading to
\begin{align*}
 \int_{\Ss_+^{N-2}} \xi \otimes \xi \diff \xi=\int_{\Ss_+^{N-2}}
 \begin{pmatrix}
 \xi_1^2 & 0 &  \\
  0 &    \ddots & 0 \\
       &  0 & \xi_{N-1}^2
 \end{pmatrix}
\diff \xi= \frac{\mathcal{H}^{N-2}(\Ss^{N-2})}{2(N-1)} Id.
\end{align*}
This concludes the proof.
\qed

\bigskip

\proof{ of Theorem \ref{pro:level}} It is a consequence of Lemma \ref{lem:second_form_propid} and a classical result in differential geometry 
for totally umbilical hypersurfaces (see e.g. \cite{Hicks} Ch. 2, page 36).
\qed

\bigskip

We have the following consequence of Lemma \ref{lem:derivative} and Theorem \ref{thm:main} (whose proof is independent of this Section):  

\begin{cor}
\label{cor:vortex}
Under the hypothesis of Theorem \ref{pro:level}, there exists a neighborhood $\om$ of $\Se$ and a diffeomorphism $t\to F(t)$ such that either
$\psi=F(|x-P|)$ for every $x\in \om$ for a point $P\in \R^N$, or $\psi=F(x\cdot \xi)$ for every $x\in \om$ for a vector $\xi \in \Ss^{N-1}$.
\end{cor}


\section{Several properties on the set of Lebesgue points}

Let $\Omega\subset \RR^N$ be an open set and $u \in L^{1}_{loc}(\Omega, \RR^N)$. 
Recall that $x_0\in \Omega$ is a \textit{Lebesgue point} of $u$ if there exists a vector $u_0\in \RR^N$ such that:
\begin{equation}\label{eq:pt_Leb}
\lim\limits_{r \to 0} \int_{B_r(x_0)}\hspace{-1.1cm} - \hspace{0.8cm}  |u(x)-u_0| \, dx= 0.
\end{equation}
In this case, we write $u(x_0)=u_0$ which is the limit of the average $\int \hspace{-0.3cm} - $ of $u$ on the ball $B_r(x_0)$ as $r\to 0$. We denote by $Leb\subset \Omega$ the set of Lebesgue points of $u$. It is well known that 
$\h^{N}(\Omega\setminus Leb)=0$ and one can replace the ball $B_r(x_0)$ by the cube $x_0+(-r,r)^N$ in the definition \eqref{eq:pt_Leb} to recover the
same set of Lebesgue points.

\proof{ of Proposition \ref{pro:curl}}
We start by proving \eqref{eq:averN} for a fixed vector $u(x)\in \Ss^{N-1}$. By carrying out a rotation if necessary,
we may assume that $u(x)=e_N$. 
Then we compute
\[
 \int_{\Ss^{N-1}}\xi \chi(x,\xi) \, d \h^{N-1} (\xi)=\int_{\Ss^{N-1} \cap \left\{ \xi_N >0 \right\} }\xi \, d \h^{N-1}(\xi)=
 \bigg(\int_{\Ss^{N-1} \cap \left\{ \xi_N >0 \right\} }\xi_N \, d \h^{N-1} (\xi)\bigg) e_N
\]
because the integrand is odd in the variables $\xi_j$ for $j=1, \dots, N-1$.
Denoting by $\xi':=(\xi_1,\dots,\xi_{N-1})$, the half sphere $\Ss^{N-1} \cap \left\{ \xi_N >0 \right\}$ is the graph of the map $\xi'\in B^{N-1}\mapsto \xi_N=\sqrt{1-|\xi'|^2}$ so that we have:
$$
\int_{\Ss^{N-1} \cap \left\{ \xi_N >0 \right\} }\xi_N \, d \h^{N-1} (\xi) 
= \int_{B^{N-1}} \sqrt{1-|\xi'|^2} \frac{d\xi'}{\sqrt{1-|\xi'|^2}}=\mathcal{H}^{N-1}(B^{N-1})=V_{N-1}. 
$$
The second statement naturally reduces (by a slicing argument) to the case of dimension $N=2$. In that case, for any $\varphi \in C^\infty_c(\Omega)$, we have $\nabla\times u=\partial_1 u_2-\partial_2 u_1$ and 
\begin{align*}
\int_{\Omega}\varphi \nabla \times u \diff x=-\int_{\Omega}\nabla^\perp \varphi \cdot u  \diff x &\stackrel{\eqref{eq:aver_form}}{=}
\frac1 2\int_{\Omega}  \int_{\Ss^1}\nabla \varphi \cdot \xi^{\perp} \chi(x,\xi)   \diff \mathcal{H}^1(\xi) \diff x \\
&=\frac1 2 \int_{\Ss^1}  \diff \mathcal{H}^1(\xi) \int_{\Omega} \nabla \varphi \cdot \xi^\perp \chi(x,\xi) \diff x \stackrel{\eqref{eq:kinet2D}}{=}0.
\end{align*}
\qed

The following lemma yields the relation between the Lebesgue points of $u$ and Lebesgue points of the functions $\{\chi(.,\xi)\}_{\xi\in \Ss^{N-1}}$ defined 
at \eqref{eq:weak_cara}.

\begin{lem}\label{lem:lebesgue}
Let $\Omega\subset \RR^N$ be an open set and $u\in L^1_{loc}(\Omega, \R^{N})$. 
\begin{enumerate}
\item[(i)] If $|u|=1$ a.e. in $\Omega$ and $x_0$ is a Lebesgue point of $\chi(.,\xi)$ for almost every $\xi\in \Ss^{N-1}$, then $x_0$ is a Lebesgue point of $u$ and \eqref{eq:averN} holds at $x_0$.

\item[(ii)] Let $x_0$ be a Lebesgue point of $u$ and $\xi\in \Ss^{N-1}$. If $u(x_0) \cdot \xi \neq 0$, then $x_0$ is a Lebesgue point of 
$\chi(.,\xi)$. Conversely, if
$x_0$ is a Lebesgue point of $\chi(.,\xi)$ with $\chi(x_0, \xi)=1$ (resp. $=0$), then { $u(x_0) \cdot \xi \geq 0$ (resp. $\leq 0$). }
\end{enumerate}

\end{lem}
\proof{} For proving $(i)$, we apply Proposition \ref{pro:curl}. Indeed, if $x_0$ is a Lebesgue point of ${\chi}(.,\xi)$ for a.e. $\xi\in \Ss^{N-1}$, then Fubini's theorem implies:
\begin{align*}
&\int_{B_r(x_0)}\hspace{-1.1cm} - \hspace{0.8cm} \bigg|u(x) - \frac 1 {V_{N-1}}  \int_{\Ss^{N-1}} \xi  {\chi}(x_0,\xi) \, d \h^{N-1}(\xi)\bigg| \diff x \\
\stackrel{\eqref{eq:averN}}{\leq} & \frac 1 {V_{N-1}} \int_{B_r(x_0)}\hspace{-1.1cm} - \hspace{0.8cm} \int_{\Ss^{N-1}} \bigg| \xi \big({\chi}(x,\xi)-{\chi}(x_0,\xi)\big) \bigg| \diff \h^{N-1}(\xi) \diff x \\
\leq & \frac 1 {V_{N-1}} \int_{\Ss^{N-1}} \bigg( \int_{B_r(x_0)}\hspace{-1.1cm} - \hspace{0.8cm} | {\chi}(x,\xi)-{\chi}(x_0,\xi)| \diff x \bigg)\diff \h^{N-1}(\xi)\stackrel{r\to 0}{\longrightarrow}0, 
\end{align*}
where we used the dominated convergence theorem. 

Next we prove $(ii)$. We treat the case $u(x_0)\cdot \xi>0$. For that, we have:
\begin{align*}
\int_{B_r(x_0)}\big| {\chi}(x,\xi)-1\big| \diff x & = 
\frac{1}{u(x_0)\cdot \xi} \int_{B_r(x_0)\cap \left\{u\cdot \xi\leq  0  \right\}} u(x_0)\cdot \xi \diff x \\
 & \leq \frac{1}{u(x_0)\cdot \xi}  \int_{B_r(x_0)\cap \left\{u\cdot \xi\leq  0  \right\}} \underbrace{(u(x_0)\cdot \xi-u(x)\cdot \xi)}_{\geq u(x_0)\cdot \xi>0} \diff x\\
& \leq \frac{1}{u(x_0)\cdot \xi}  \int_{B_r(x_0)}|u(x)-u(x_0)| \diff x.
\end{align*}
Since $x_0$ is a Lebesgue point of $u$, it follows that $x_0$ is a Lebesgue point for ${\chi}(\cdot,\xi)$ with ${\chi}(x_0,\xi)=1$. The case $u(x_0)\cdot \xi<0$ can be shown similarly and obtain ${\chi}(x_0,\xi)=0$. The last statement is a direct consequence of the above lines (using a contradiction argument).
\qed

\begin{rem} $a)$ Note that the condition $u(x_0) \cdot \xi \neq 0$ is essential in Lemma \ref{lem:lebesgue} $(ii)$. Indeed, if one considers the vortex vector field 
$u(x)=\frac{x}{|x|}$ for $x\in \R^N\setminus \{0\}$, then for every $\xi \in \Ss^{N-1}$, any point $x_0 \in \xi^\perp\setminus\{0\}$ is a Lebesgue point of $u$ (because $u$ is smooth around $x_0$) and satisfies $$u(x_0)\cdot \xi=0,$$ but $x_0$ 
is not a Lebesgue point of $\chi(\cdot, \xi)$ because
\[
 \int_{B_r(x_0)}\hspace{-1.1cm} - \hspace{0.8cm} \bigg|{\chi}(x,\xi)-  \int_{B_r(x_0)}\hspace{-1.1cm} - \hspace{0.8cm} {\chi}(\cdot,\xi)\bigg|\diff x=\int_{B_r(x_0)}\hspace{-1.1cm} - \hspace{0.8cm}\frac 1 2\diff x \nrightarrow 0 \quad \textrm{ as } r\to 0,
\]
where we used that
\[
\int_{B_r(x_0)}\hspace{-1.1cm} - \hspace{0.8cm} {\chi}(x,\xi) \diff x= \frac{\h^{N}(\left\{x\in B_r(x_0)\, :\, x\cdot \xi >0 \right\})}{\h^{N}(B_r(x_0))} \stackrel{x=y+x_0}{=}
\frac{\h^{N}(\left\{y\in B_r(0)\, :\, y\cdot \xi >0 \right\})}{\h^{N}(B_r(0))} = \frac{1}{2}.
\]
{ $b)$ Note that in Lemma \ref{lem:lebesgue} $(ii)$, one cannot conclude in general that $u(x_0) \cdot \xi > 0$ provided that ${\chi}(x_0,\xi)=1$. Indeed, consider for example $\xi=e_N$, $u(x)\cdot \xi=u_N(x):=|x|$ for $x\in \R^N$ and set $x_0=0$; then $\chi(\cdot, \xi)=1$ in $\R^N\setminus\{x_0\}$, $x_0$ is a Lebesgue point of $u_N$ and $\chi(\cdot, \xi)$ with $\chi(x_0, \xi)=1$, but $u_N(x_0)=0$.}
\end{rem}

We now prove one of the key tools in the proof of Theorem \ref{thm:main} that mimics the relation of ordering of level sets of a stream function when \eqref{eq:kineN}
holds true. It is a generalization of Proposition~3.1 in \cite{JOP02} to the case of dimension $N$:

\begin{pro}[Ordering]
\label{pro:preserv_prod}
Let $N\geq 2$, $\Omega\subset \RR^N$ be an open set and $u\in L^1_{loc}(\Omega, \R^{N})$ satisfying the kinetic formulation \eqref{eq:kineN}. Assume that $y, z \in Leb$ are two different Lebesgue points of $u$ such that the closed segment $[yz]$ is included in $\Omega$. Then for every direction $\xi \in \Ss^{N-1}$ with $\xi \in  (z-y)^{\perp}$, we have: 
\be \label{eq:cas_scalaire}
  u(y) \cdot \xi >0 \textrm{ (resp. $<0$)} \quad \Rightarrow \quad u(z) \cdot \xi { \geq 0} \textrm{ (resp. ${ \leq 0}$)};
\ee 
moreover, $y$ and $z$ are Lebesgue points of $\chi(\cdot, \xi)$ and $\chi(y, \xi)=\chi(z, \xi)$. 
As a consequence, if $u\neq 0$ a.e. in $\Omega$, then we have for a.e. $y\in \Omega$, $\h^{N-1}$-a.e. $\xi \in \Ss^{N-1}$ and $\h^{N-1}$-a.e. $v\in \xi^\perp$ with the segment $[y, y+v]\subset \Omega$ that $y$ and $y+v$ are Lebesgue points of $u$ and
\be \label{eq:egal_chi}
\chi(y, \xi)=\chi(y+v, \xi).
\ee
\end{pro}
\proof{}
First, we consider the case $u(y) \cdot \xi > 0$. By Lemma~\ref{lem:lebesgue} $(ii)$,
$y$ is a Lebesgue point of ${\chi}(.,\xi)$ and ${\chi}(y,\xi)=1$.
Let $$\{\rho_\eps(\cdot)=\frac{1}{\eps^N}\rho(\frac{\cdot}{\eps})\}_{\eps>0}$$ be a standard family of mollifiers where $\rho$ is a nonnegative radial smooth function having as support the unit ball $\supp \rho=B_1\subset \R^N$ and $\int_{B_1} \rho\, dx=1$. Set the convoluted function
\[
{\chi}_{\varepsilon}:= \rho_{\varepsilon} *{\chi}(.,\xi)
\]
in a neighborhood $\omega\subset \Omega$ of the segment $[yz]$ for $\eps>0$ sufficiently small. Then ${\chi}_{\varepsilon}$ is smooth in $\omega$ and for every Lebesgue point 
$x\in \omega$ of ${\chi}(.,\xi)$ we have ${\chi}_{\varepsilon}(x) \to {\chi}(x,\xi)$ as $\eps\to 0$ because
\begin{align*}
|{\chi}_{\varepsilon}(x)-{\chi}(x,\xi)|&=\bigg|\int_{B_\varepsilon(0)}\big({\chi}(x-\tilde{x},\xi)-{\chi}(x,\xi)\big)\rho_\eps(\tilde x) \diff \tilde x\bigg|\\
& \leq \frac{\sup \rho}{\eps^N}\int_{B_\varepsilon(0)}|{\chi}(x-\tilde{x},\xi)-{\chi}(x,\xi)| \diff \tilde{x}\leq C \int_{B_\eps(x)}\hspace{-1cm} - \hspace{0.8cm} |{\chi}(\tilde{y},\xi)-{\chi}(x,\xi)| \diff \tilde{y}\stackrel{\eps\to 0}{\to} 0.
\end{align*}
 In particular, $\lim_{\eps\to 0}\chi_\eps(y)=\chi(y, \xi)=1$.
Let $v=z-y$. We show that $\chi(y+v,\xi)=1$. For that, we have $v\in \xi^\perp$ and 
$$v\cdot \nabla_x \chi_\eps=v\cdot \nabla_x  {\chi}(\cdot,\xi)*\rho_\eps\stackrel{\eqref{eq:kineN}}{=}0\quad \textrm{in}\quad \omega.$$
Then
\[
\chi_\eps(y+v)-\chi_\eps(y)=\int_0^1 v \cdot \nabla_x \chi_\eps(y+tv)\, dt=0
\]
so that $\lim_{\eps\to 0}\chi_\eps(z)=\lim_{\eps\to 0}\chi_\eps(y)=\chi(y, \xi)=1$. This implies that $u(z)\cdot \xi\geq 0$. Assume by contradiction that $u(z)\cdot \xi<0$. By Lemma \ref{lem:lebesgue} ii), $z$ is a Lebesgue point of ${\chi}(.,\xi)$ with $\chi(z, \xi)=0$ so that 
$\lim_{\eps\to 0}\chi_\eps(z)=\chi(z, \xi)=0$ which contradicts the above statement. 
We prove now the following:

\medskip

{
\nd {\bf Claim:} If ${\chi}_{\varepsilon}(z) \to 1$ as $\eps\to 0$,  then $z$ is a Lebesgue point of ${\chi}(.,\xi)$ with ${\chi}(z,\xi)=1$. }
\proof{ of Claim} Let $\{\eps_k\}$ be a sequence converging to $0$ as $k\to \infty$. For $k$ large enough, we define $f_k:B_1\to \{0, 1\}$ by 
$f_k(x)=\chi(z-\eps_k x, \xi)$ for every $x\in B_1$. Then the sequence $\{f_k\}$ is bounded in $L^2(B_1)$ and up to a subsequence, $f_k\rightharpoonup f$ weakly in $L^2(B_1)$ where the limit $f:B_1\to \R$ has the range inside $[0,1]$. Therefore, we have for our smooth mollifier $\rho\in L^2(B_1)$ that
$$\int_{B_1} \rho\, f_k\, dx\to \int_{B_1} \rho\, f\, dx \quad \textrm{as} \quad k\to \infty.$$ 
Note now that by the change of variable $\tilde{x}=z-\eps_k x$ we obtain by our assumption:
$$\int_{B_1} \rho(x) f_k(x)\, dx=\int_{B_{\eps_k}(z)}\rho_{\eps_k}(z-\tilde x) \chi(\tilde x, \xi) \, d\tilde x={\chi}_{\varepsilon_k}(z)\to 1 \quad \textrm{as} \quad k\to \infty,$$
therefore, $\int_{B_1} \rho\, f\, dx=1$. Since $1$ is the maximal value of $f$ and $\rho$ is nonnegative with the integral on $B_1$ equal to $1$, we deduce that $f=1$ in $\supp\rho=B_1$. It follows by changing the variable $\tilde{x}=z-\eps_k x$:
$$\int_{B_{\eps_k}(z)}\hspace{-1.1cm} - \hspace{0.8cm} |{\chi}(\tilde x,\xi)-1| \diff \tilde x=1-\int_{B_1(0)}\hspace{-0.95cm} - \hspace{0.8cm} f_k(x)\, dx\to 0
\quad \textrm{as} \quad k\to \infty,$$
because $f_k\rightharpoonup 1$ weakly in $L^2(B_1)$. Since the limit is unique for every subsequence $\eps_k\to 0$, we conclude that $z$ is a Lebesgue point of ${\chi}(.,\xi)$ with ${\chi}(z,\xi)=1$ which proves the claim.\\

\nd For the case $u(y) \cdot \xi < 0$ (i.e., ${\chi}(y,\xi)=0$ by Lemma \ref{lem:lebesgue} $(ii)$), one applies the above argument by replacing $\xi$ by $-\xi$ and obtain that  $z$ is a Lebesgue point of ${\chi}(.,-\xi)$ with ${\chi}(z,-\xi)=1$. It follows that $z$ is a Lebesgue point of ${\chi}(.,\xi)$ with ${\chi}(z,\xi)=0$ because
$$\int_{B_r(z)}\hspace{-0.95cm} - \hspace{0.8cm} |{\chi}(x,\xi)| \diff x\leq \frac{\h^{N}(\left\{x\in B_r(z)\, :\, u(x)\cdot \xi \geq 0 \right\})}{\h^{N}(B_r(z))}=1-\int_{B_r(z)}\hspace{-0.95cm} - \hspace{0.8cm} {\chi}(x,-\xi) \diff x\to 0$$
as $r\to 0$. One also concludes that $u(z)\cdot \xi\leq 0$ by Lemma \ref{lem:lebesgue} $(ii)$.

For the last statement, we have for a.e. $y\in \Omega$, $y$ is a Lebesque point of $u$ with $u(y)\neq 0$. Then for $\h^{N-1}$-a.e. direction $\xi \in \Ss^{N-1}$, we have that $u(y)\cdot \xi\neq 0$ and $y+v$ is a Lebesgue point of $u$ for $\h^{N-1}$-a.e. $v\in \xi^\perp$ with the segment $[y, y+v]\subset \Omega$. By the above argument, we conclude to \eqref{eq:egal_chi}.
\qed

\section{Notion of trace on lines}


The $H^{1/2}$-regularity for $N$-dimensional unit length vector fields $u$ satisfying the kinetic formulation \eqref{eq:kineN} (see \cite{Golse}) is a priori not enough to define the notion of trace of $u$ on one-dimensional lines. However, 
using the ideas in \cite{JOP02} for dimension $2$, we will define a notion of trace of $u$ on segments (in the sense of Lebesgue points) in any dimension $N\geq 2$. 

\begin{pro}[Trace]\label{pro:trace}
Let $N\geq 2$, $\Omega\subset \RR^N$ be an open set and $u:\Omega\to \Ss^{N-1}$ be a Lebesgue measurable vector field satisfying the kinetic formulation \eqref{eq:kineN}. Assume that the segment $$L:=\left\{0 \right\}^{N-1}\times [-1,1] \textrm{ is included in $\Omega$}.$$ Then there exists a Lebesgue measurable function  
$\tilde{u}:(-1,1)\to \R^N$ such that 
\be \label{eq:trace_form}
\lim\limits_{r \to 0}  \int_{(-r,r)^{N-1}} \hspace{-1.55cm} - \hspace{1.2cm}  \int_{-1}^{1}
\abs{u(x', x_N)-\tilde{u}(x_N)} \diff x'\diff x_N= 0,
\ee
where $x=(x', x_N)$, $x'=(x_1, \dots, x_{N-1})$. 
Moreover, for $\mathcal{H}^1$-a.e. $x_N \in (-1,1)$,
\begin{align}\label{eq:trace_droite}
\tilde{u}(x_N)= \lim\limits_{r \to 0} \int_{(-r,r)^{N-1}} \hspace{-1.55cm} - \hspace{1.2cm}  u(x', x_N) \, dx' \quad \text{and} \quad |\tilde{u}(x_N)|=1.
\end{align}
\footnote{$Leb$ is the set of Lebesgue points of $u$ in $\Omega$.}Finally, every Lebesgue point $x \in Leb$ of $u$ lying inside $L$ is a Lebesgue point of $\tilde u$ and
$u(x)=\tilde{u}(x_N)$. The vector field $\tilde u$ is called the trace of $u$ on the segment $L$.
\end{pro}
\proof{} To simplify the writing, we assume that $\Omega=\RR^N$. We divide the proof in several steps:

\medskip

\nd {\it Step 1: Defining the one-dimensional function $\tilde\chi(\cdot, \xi)$ for suitable directions $\xi\in \Ss^{N-1}$}. Let $\cal D$ be the set of directions $\xi \in \Ss^{N-1}$ such that $\xi_N \neq 0$ and \eqref{eq:egal_chi} holds true for the triple $(y, y+v, \xi)$ for a.e. $y\in \Omega$ and $\h^{N-1}$-a.e. $v\in \xi^\perp$ (with the segment $[y, y+v]\subset \Omega$ where $y$ and $y+v$ are Lebesgue points of $u$). By Proposition \ref{pro:preserv_prod}, we know that $\cal D$ covers $\Ss^{N-1}$ up to a set of zero $\h^{N-1}$-measure. For such a direction $\xi\in {\cal D}$, we can choose a point $y_\xi\in \Omega$ (in a neighborhood of $L$) such that the map $\xi\in {\cal D}\mapsto y_\xi\in \Om$ is Lebesgue measurable, the point $y_\xi+t\xi\in \Omega$ is a Lebesgue point of $\chi(\cdot,\xi)$ for $\h^1$-a.e. $t\in \R$, the function $t\mapsto \chi(y_\xi+t \xi,\xi)$
is $\h^1$-measurable (by Fubini's theorem) and \eqref{eq:egal_chi} holds true for the triple $(y_\xi+t\xi, y_\xi+t\xi+v, \xi)$ for $\h^{N-1}$-a.e. $v\in \xi^\perp$ and 
$\h^{1}$-a.e. $t$. Set the one-dimensional function
$$s\mapsto \tilde{\chi}(s, \xi):=\chi\big(y_\xi+(s-y_\xi \cdot \xi) \xi,\xi\big) \in \{0,1\}.$$
Then we have that for a.e. $x\in\Omega$ in a neighborhood of $L$:
\be \label{eq:12}
\tilde{\chi}(x \cdot \xi,\xi) = \chi(y_\xi-y_\xi\cdot \xi \xi+x\cdot \xi \xi,\xi)\stackrel{\eqref{eq:egal_chi}}{=}{\chi}(x,\xi),
\ee
because $v=y_\xi-y_\xi\cdot \xi \xi+x\cdot \xi \xi-x\in \xi^\perp$.

\medskip

\nd {\it Step 2: For $\xi\in {\cal D}$ and for every Lebesgue point $P=(0, \dots, 0, P_N)\in L$ of ${\chi}(\cdot,\xi)$ with $P_N\in (-1,1)$, the point $P\cdot \xi$ is a Lebesgue point of $\tilde\chi(\cdot, \xi)$ and $\tilde\chi(P_N\xi_N, \xi)={\chi}(P,\xi)$.}  Indeed, since $\xi_N\neq 0$, we have:
\begin{align*}
 \int_{P_N\xi_N-r|\xi_N|}^{P_N\xi_N+r|\xi_N|}\hspace{-1.99cm} - \hspace{1.63cm}\big|\tilde\chi(t, \xi) &-{\chi}(P,\xi)  \big| \diff t 
\stackrel{t=\tilde{x}_N\xi_N}{=}  \int_{(-r,r)^{N-1}} \hspace{-1.55cm} - \hspace{1.2cm}\diff x'  \int_{P_N-r}^{P_N+r} \hspace{-1.09cm} - \hspace{0.75cm}  \big|\tilde\chi
(\tilde{x}_N\xi_N, \xi)-{\chi}(P,\xi)  \big| \diff \tilde{x}_N \\
  \stackrel{x'\cdot \xi'+{x}_N\xi_N=\tilde{x}_N\xi_N}{=} & \int_{(-r,r)^{N-1}} \hspace{-1.55cm} - \hspace{1.2cm}\diff x' \int_{P_N-\frac{x'\cdot \xi'}{\xi_N}-r}^{P_N-\frac{x'\cdot \xi'}{\xi_N}+r} \hspace{-1.92cm} - \hspace{1.63cm}  \big|\underbrace{\tilde \chi(x'\cdot \xi'+{x}_N\xi_N, \xi)}_{\stackrel{\eqref{eq:12}}{=} {\chi}(x,\xi)}-{\chi}(P,\xi)  \big| \diff {x}_N\\
\hspace{-3cm}   \leq &  \int_{(-r,r)^{N-1}} \hspace{-1.55cm} - \hspace{1.2cm}\diff x'\frac{1}{2r} \int_{P_N-\rti}^{P_N+\rti}\big|{\chi}(x,\xi)-{\chi}(P,\xi)  \big|\, dx_N\\
\hspace{-3cm}   \leq &C \int_{P+(-\rti, \rti)^{N}} \hspace{-1.65cm} - \hspace{1.4cm}\big|{\chi}(x,\xi)-{\chi}(P,\xi)  \big|\, dx\, \, 
    {\to} \, \, 0 
  \quad \textrm{as $r\to 0$}
\end{align*}
where we used that $|x'\cdot \xi'|\leq r\sqrt{N-1}$ for $x'\in (-r,r)^{N-1}$ and $\rti=\big(\frac{\sqrt{N-1}}{|\xi_N|}+1\big)r$. Thus, $P_N\xi_N$ is a Lebesgue point of 
$\tilde\chi(\cdot, \xi)$. In particular, we have by Fubini's theorem for every $\alpha>0$:
\begin{align}
\nonumber
& \int_{-\alpha r}^{\alpha r} \hspace{-0.8cm} - \hspace{0.5cm} \, d\tit \int_{P_N\xi_N-r|\xi_N|+\tit}^{P_N\xi_N+r|\xi_N|+\tit}\hspace{-2.3cm} - \hspace{2cm}\big|\tilde\chi(t, \xi)-\tilde\chi(P_N\xi_N, \xi) \big| \diff t \\
\nonumber
&=\frac{1}{4\alpha |\xi_N| r^2}  \int_{-\alpha r}^{\alpha r} \int_{P_N\xi_N-r(|\xi_N|+\alpha)}^{P_N\xi_N+r(|\xi_N|+\alpha)} \big|\tilde\chi(t, \xi)-\tilde\chi(P_N\xi_N, \xi) \big| 
\, {\mathds{1}}_{(P_N\xi_N-r|\xi_N|+\tit, P_N\xi_N+r|\xi_N|+\tit)}(t)\, dt\, d\tit\\
\nonumber
&=\frac{1}{4\alpha |\xi_N| r^2}   \int_{P_N\xi_N-r(|\xi_N|+\alpha)}^{P_N\xi_N+r(|\xi_N|+\alpha)} \big|\tilde\chi(t, \xi)-\tilde\chi(P_N\xi_N, \xi) \big| \, dt
\int_{-\alpha r}^{\alpha r} {\mathds{1}}_{(-P_N\xi_N-r|\xi_N|+t, -P_N\xi_N+r|\xi_N|+t)}(\tit)\, d\tit \\
\label{eq:need}
&\leq \frac{1}{2|\xi_N| r}   \int_{P_N\xi_N-r(|\xi_N|+\alpha)}^{P_N\xi_N+r(|\xi_N|+\alpha)} \big|\tilde\chi(t, \xi)-\tilde\chi(P_N\xi_N, \xi) \big| \, dt \, \, {\to} \, \, 0 
  \quad \textrm{as $r\to 0$}.
\end{align}

\medskip

\nd {\it Step 3: Proof of \eqref{eq:trace_form}}. 
For $\xi\in {\cal D}$, we have for small $r>0$:
\begin{align*}
\int_{(-r,r)^{N-1}} \hspace{-1.55cm} - \hspace{1.2cm} \int_{-1}^{1} 
&\abs{{\chi}(x,\xi)-\tilde{{\chi}}(x_N  \xi_N,\xi)}  \diff x' \diff x_N\\  
   \stackrel{\eqref{eq:12}}{=} & \int_{(-r,r)^{N-1}} \hspace{-1.55cm} - \hspace{1.2cm} \int_{-1}^{1}
\abs{\tilde{{\chi}}(x'\cdot \xi'+x_N\xi_N,\xi)-\tilde{{\chi}}(x_N \xi_N,\xi)}  \diff x' \diff x_N \\
  \stackrel{t=x_N\xi_N}{\leq} & \frac{1}{|\xi_N|}\underset{\abs{\tilde t} \leq r\sqrt{N-1}}{\sup} \int_{-|\xi_N|}^{|\xi_N|} 
\abs{\tilde{{\chi}}(t+\tilde t,\xi)-\tilde{{\chi}}(t,\xi)} \diff t
\end{align*}
because $|x'\cdot \xi'|\leq r\sqrt{N-1}$. 
Since the one-dimensional function $t \mapsto \tilde{{\chi}}(t ,\xi)$ belongs to $L^\infty$, its $L^1$-modulus of continuity present in the above RHS tends to $0$ as $r\to 0$ which leads to the following:
\begin{align*}
\lim\limits_{r \to 0} \int_{(-r,r)^{N-1}} \hspace{-1.55cm} - \hspace{1.2cm} \int_{-1}^{1}
\abs{{\chi}(x,\xi)-\tilde{{\chi}}(x_N \xi_N,\xi)} \diff x'  \diff x_N =0.
\end{align*}
This formula can be interpreted as the notion of trace of $\chi(\cdot, \xi)$ on the segment $L$ and yields \eqref{eq:trace_form}. Indeed, due to \eqref{eq:averN}, 
we set for a.e. $x_N\in (-1, 1)$:
\[
 \tilde{u}(x_N)=\frac 1 {V_{N-1}} \int_{\Ss^{N-1}} \xi \tilde{{\chi}}(x_N \xi_N,\xi) \diff \h^{N-1} (\xi)
\]
and we obtain by Fubini's theorem:
\begin{align*}
\int_{(-r,r)^{N-1}} \hspace{-1.55cm} - \hspace{1.2cm} & \int_{-1}^{1}
\abs{u(x', x_N)-\tilde{u}(x_N)} \diff x' \diff x_N \\
\stackrel{\eqref{eq:averN}}{\leq} & \frac 1 {V_{N-1}} \int_{\Ss^{N-1}} \bigg(\int_{(-r,r)^{N-1}} \hspace{-1.55cm} - \hspace{1.2cm} \int_{-1}^{1}
\abs{{\chi}(x,\xi)-\tilde{{\chi}}(x_N\xi_N,\xi)} \diff x' \diff x_N \bigg) \diff \h^{N-1} (\xi)\to 0 \quad \textrm{ as } r\to 0,
\end{align*}
where we used the dominated convergence theorem. 

\medskip

\nd {\it Step 4: Proof of \eqref{eq:trace_droite}}. By Step 3, we deduce that:
\[
\int_{(-r,r)^{N-1}} \hspace{-1.55cm} - \hspace{1.2cm}  u(x', \cdot) \diff x' \xrightarrow[r \to 0]{}  \tilde{u} \quad \text{in} \; L^{1}((-1,1));
\]
therefore, the first statement in \eqref{eq:trace_droite} follows immediately.
Moreover, 
\begin{align*}
 \int_{-1}^1\bigg| |\tilde{u}(x_N)|-1 \bigg| \diff x_N 
 =  & \int_{(-r,r)^{N-1}} \hspace{-1.55cm} - \hspace{1.2cm} \int_{-1}^1 \bigg| |\tilde{u}(x_N)|-|u(x', x_N)| \bigg| \diff x' \diff x_N  \\
  \leq & \int_{(-r,r)^{N-1}} \hspace{-1.55cm} - \hspace{1.2cm}\int_{-1}^{1} \abs{\tilde{u}(x_N)-u(x', x_N)}  \diff x' \diff x_N\stackrel{\eqref{eq:trace_form}}{\to} 0 
  \quad \textrm{as $r\to 0$;}
\end{align*}
thus, $|\tilde{u}(x_N)|=1$ for $\mathcal{H}^1$-a.e. $x_N\in (-1,1)$.

\medskip

\nd {\it Step 5: Conclusion}. Let $P=(0, \dots, 0, P_N) \in Leb$ be a Lebesgue point of $u$ with $P_N\in (-1, 1)$. We want to show that $P_N$ is a Lebesgue point of $\tilde u$ and $\tilde u(P_N)=u(P)$. For that, we know by Lemma \ref{lem:lebesgue} that $P$ is a Lebesgue point of $\chi(\cdot, \xi)$ for every direction $\xi\in \Ss^{N-1}$ with $u(P)\cdot \xi\neq 0$. If in addition $\xi\in \cal D$, we know by Step 2 that $P\cdot \xi$ is also a Lebesgue point of $\tilde\chi(\cdot, \xi)$. By the same argument as at Step 3, we have:
\begin{align*}
\int_{P+(-r,r)^{N}} \hspace{-1.63cm} - \hspace{1.3cm} &
\abs{u(x', x_N)-\tilde{u}(x_N)}  \diff x' \diff x_N \\
\leq & \frac 1 {V_{N-1}} \int_{\Ss^{N-1}} \int_{P+(-r,r)^{N}} \hspace{-1.63cm} - \hspace{1.3cm} \bigg|\underbrace{{\chi}(x,\xi)}_{\stackrel{\eqref{eq:12}}{=} \tilde{{\chi}}(x'\cdot \xi'+x_N\xi_N,\xi)}-\tilde{{\chi}}(x_N\xi_N,\xi)\bigg|  \diff x' \diff x_N \diff \h^{N-1} (\xi)\\
\leq & \frac 1 {V_{N-1}} \int_{\Ss^{N-1}}  \diff \h^{N-1} (\xi) \bigg[\int_{P+(-r,r)^{N}} \hspace{-1.63cm} - \hspace{1.3cm} \bigg|\tilde{\chi}(x'\cdot \xi'+x_N\xi_N,\xi)-\tilde{\chi}(P_N\xi_N,\xi)\bigg|\, dx \\
&\hspace{2cm} + \int_{P_N-r}^{P_N+r} \hspace{-1.15cm} - \hspace{0.95cm} \bigg| \tilde{{\chi}}(x_N\xi_N,\xi)-\tilde{\chi}(P_N\xi_N,\xi)\bigg| \diff x_N\bigg]\\
\leq & \frac 1 {V_{N-1}} \int_{\Ss^{N-1}} \, d\h^{N-1} (\xi) \int_{(-r,r)^{N-1}} \hspace{-1.53cm} - \hspace{1.3cm} \, dx' \int_{P_N\xi_N-r|\xi_N|+x'\cdot \xi'}^{P_N\xi_N+r|\xi_N|+x'\cdot \xi'} \hspace{-2.75cm} - \hspace{2.5cm} \bigg|\tilde{\chi}(t,\xi)-\tilde{\chi}(P_N\xi_N,\xi)\bigg|\, dt \\
&\hspace{2cm} +\frac 1 {V_{N-1}} \int_{\Ss^{N-1}} \, d\h^{N-1} (\xi) \int_{P_N\xi_N-r|\xi_N|}^{P_N\xi_N+r|\xi_N|} \hspace{-1.92cm} - \hspace{1.63cm} \bigg| \tilde{{\chi}}(t,\xi)-\tilde{\chi}(P_N\xi_N,\xi)\bigg| \diff t. 
\end{align*}
Using twice the dominated convergence theorem, we conclude that the above RHS vanishes as $r\to 0$. Indeed, the second integrand converges to $0$ as $r\to 0$ by Step 2 for a.e. $\xi \in \Ss^{N-1}$. For the first integrand, we proceed as follows: for $\h^{N-1}$-a.e. direction $\xi$, we may assume that $|\xi'|> 0$ and $\xi_N\neq 0$ so that there exists a rotation 
$R'\in SO(N-1)$ with $R'\xi'=|\xi'|e_1$ and we have by the change of variable $\tilde x'=R'x'$ and $\hat{r}=r\sqrt{N-1}$:
\begin{align*}
&\int_{(-r,r)^{N-1}} \hspace{-1.58cm} - \hspace{1.3cm} \, dx' \int_{P_N\xi_N-r|\xi_N|+x'\cdot \xi'}^{P_N\xi_N+r|\xi_N|+x'\cdot \xi'} \hspace{-2.75cm} - \hspace{2.5cm} \bigg|\tilde{\chi}(t,\xi)-\tilde{\chi}(P_N\xi_N,\xi)\bigg|\, dt\\
&\leq C \int_{\{|\tilde x'|<\hat r\}} \hspace{-1.3cm} - \hspace{1.1cm} \, d\tilde x' \int_{P_N\xi_N-r|\xi_N|+\tilde x_1|\xi'|}^{P_N\xi_N+r|\xi_N|+\tilde x_1|\xi'|} \hspace{-2.9cm} - \hspace{2.6cm} \bigg|\tilde{\chi}(t,\xi)-\tilde{\chi}(P_N\xi_N,\xi)\bigg|\, dt\\
&\leq C \int_{-|\xi'| \hat r}^{|\xi'| \hat r} \hspace{-1cm} - \hspace{0.8cm} \int_{P_N\xi_N-r|\xi_N|+\tit}^{P_N\xi_N+r|\xi_N|+\tit}\hspace{-2.3cm} - \hspace{1.93cm}\big|\tilde\chi(t, \xi)-\tilde\chi(P_N\xi_N, \xi) \big| \diff t \, d\tit \, \, 
\stackrel{\eqref{eq:need}}{\to} 0 \, \,  \textrm{ as } r\to 0.
\end{align*}
\qed

\section{Proof of Theorem \ref{thm:main}}

We start by showing some preliminary results that reveal the geometric consequences of the kinetic formulation \eqref{eq:kineN}. 
The following lemma is the first step for proving that $u$ is constant along the characteristics and is reminiscent of the ideas presented in \cite{JOP02}:

\begin{lem}\label{lem:const_droite}
Let $\Omega\subset \R^N$ be an open set such that $L= \left\{0 \right\}^{N-1}\times [-1,1] \subset \Omega$ and 
$u:\Omega\to \Ss^{N-1}$ be a Lebesgue measurable vector field satisfying the kinetic formulation \eqref{eq:kineN}.
Assume that the origin
$O \in \Omega \cap Leb$ is a Lebesgue point of $u$ and $u(O)=e_N$. Then we have for every Lebesgue point $x_N \in (-1,1)$ of $\tilde u$: 
$$\tilde{u}(x_N)= \pm e_N,$$ where $\tilde{u}$ is the trace of $u$ on $L$ defined at Proposition \ref{pro:trace}.
\end{lem}

\proof{} W.l.o.g. we assume that $\Omega$ is a convex open neighborhood of $L$. By Proposition \ref{pro:trace}, we know that $O$ is also a Lebesgue point of $\tilde{u}$
and $\tilde{u}(0)=e_N$; moreover, $|\tilde u|=1$ a.e. in $(-1,1)$. Let $x_N \in (-1,1)\setminus \{0\}$ be a Lebesgue point of $\tilde{u}$ such that $\h^{N-1}$-a.e. $z\in \Omega \cap \big(x_N e_N+e_N^\perp\big)$ is a Lebesgue point of $u$
and that the following holds true (see Proposition~\ref{pro:trace}):
\be
\label{eq:20}
\lim\limits_{r \to 0}  \int_{(-r,r)^{N-1}}  \hspace{-1.5cm} - \hspace{1.2cm} \big|u(x', x_N)-\tilde{u}(x_N)\big| \diff x'=0.\ee
Our goal is to prove that the component $\tilde{u}_i(x_N)$ of $\tilde{u}(x_N)$ in direction $e_i$ vanishes for every $i=1, \dots, N-1$. For that, we follow the ideas in \cite{JOP02}. 
Let $\eps>0$ be small and denote the following subsets $E_i^-$ and $E_i^+$ (depending on $\eps$) of the hyperplane $\big(x_Ne_N+e_N^\perp\big)$ for $1 \leq i\leq N-1$: 
\[
E_i^{\pm}= \left\{z \in \Omega \cap Leb \;:\; z_N= x_N, \; \varepsilon |x_N| \geq \pm z_i > 0 \right\}. 
\] 
By our assumption, these sets $E_i^{\pm}$ contain many points (e.g., for $i=1$, $E_1^{+}$ covers the $N-1$ parallelepiped $(0,r)\times (-r,r)^{N-2}\times \{x_N\}$ up to a set of zero $\h^{N-1}$-measure, for $r<\eps$). For $z \in E_i^+$, we set $y=-z_i e_N + x_N e_i$ if $x_N>0$ (resp., $y=z_i e_N -x_N e_i$ if $x_N<0$). Obviously, $z \cdot y =0$, i.e., $y\in z^\perp$. By convexity of 
$\Omega$, the segment $[Oz]\subset \Omega$ so that by Proposition \ref{pro:preserv_prod} we have if $x_N>0$ (resp. $x_N<0$), then
$u(O) \cdot y=-z_i <0$ (resp. $u(O) \cdot y=z_i >0$) so that $u(z) \cdot y \leq 0$ (resp. $\geq 0$). It follows that 
         $$u_i(z) \leq \frac{z_i}{x_N}u_N(z) \leq \varepsilon \quad (\textrm{resp.}\, u_i(z) \geq \frac{-z_i}{|x_N|}u_N(z) \geq -\varepsilon),
$$
because $|u_N(z)|\leq 1$. Similarly, for $z \in E_i^-$, one uses $y=z_i e_N -x_N e_i$ if $x_N>0$ (resp. $y=-z_i e_N +x_N e_i$ if $x_N<0$) and deduces
that  $u_i(z) \geq -\varepsilon$  if $x_N>0$ (resp. $u_i(z) \leq \varepsilon$  if $x_N<0$). We conclude that $\tilde{u}_i(x_N)\in [-\eps, \eps]$. Indeed, let us set $i=1$ for simplicity of writing; by \eqref{eq:20}, we have
  \begin{align*}
 \tilde{u}_1(x_N) = \lim\limits_{r \to 0}  \int_{ (0,r) \times (-r,r)^{N-2}}  \hspace{-2.4cm} - \hspace{2.2cm} u_1(x', x_N) \diff x'\leq \eps \textrm{ if } x_N>0 \quad (\textrm{resp. } \geq -\eps \textrm{ if } x_N<0)
\end{align*}
and also,
\begin{align*}
 \tilde{u}_1(x_N) = \lim\limits_{r \to 0} \int_{ (-r,0) \times (-r,r)^{N-2}}  \hspace{-2.6cm} - \hspace{2.4cm} u_1(x', x_N) \diff x'\geq -\eps \textrm{ if } x_N>0 \quad (\textrm{resp. } \leq \eps \textrm{ if } x_N<0).
\end{align*}
Passing to the limit $\varepsilon\to 0$, we conclude that $\tilde{u}_i(x_N)=0$ for $i=1$ (similarly, for every $1\leq i\leq N-1$). Obviously, $\h^{1}$-a.e. $x_N\in (-1,1)$ satisfies this property. As a consequence, if $P_N\in (-1,1)$ is a Lebesgue point of $\tilde{u}$ then for every $1 \leq i\leq N-1$:
$$\tilde{u}_i(P_N)=  \lim\limits_{r \to 0} \int_{P_N-r}^{P_N+r}  \hspace{-1.1cm} - \hspace{0.75cm} \tilde{u}_i(x_N) \diff x_N=0.$$
Since $|\tilde u(P_N)|=1$, we deduce that $\tilde u_N(P_N)=\pm 1$, i.e., $\tilde u(P_N)=\pm e_N$. 
\qed

\bigskip

We now prove the main result: 
\proof{ of Theorem \ref{thm:main}} We first treat the case $\Omega$ is a ball and then the general case of a connected open set.

\medskip

\nd {\bf Case I: $\Omega$ is a ball}. Since $u$ is not a constant vector field, there exist two Lebesgue points $P_0, P_1\in \Omega \cap Leb$ of $u$ such that 
$$u(P_0) \neq u(P_1).$$ Let $D_0$ (resp. $D_1$) be the line directed by $u(P_0)$ (resp. $u(P_1)$) that passes through $P_0$ (resp. $P_1$). 

\medskip

\nd {\it Step 1: We show that $D_0$ and $D_1$ are coplanar}. Assume by contradiction that $D_0$ and $D_1$ are not coplanar, in particular $|u(P_0) \cdot u(P_1)|<1$. Set $A\in D_0$ 
and $B\in D_1$ such that
\[
0< |A-B|= \min_{x \in D_0, y \in D_1} |x-y|.
\]
\begin{figure}[htbp]
\center
\includegraphics[scale=0.5,
width=0.5\textwidth]{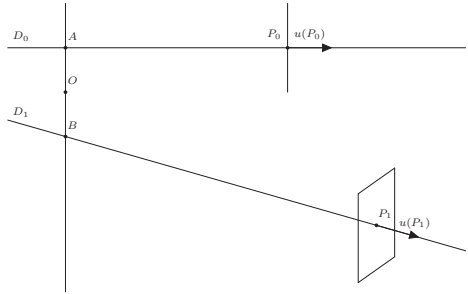} \caption{Two noncoplanar lines $D_0$ and $D_1$}
\label{no_coplan} 
\end{figure}
Obviously, the segment $[AB]$ is orthogonal to $D_0$ and $D_1$. Set $O$ be the middle point of the segment $[AB]$ (see Figure \ref{no_coplan}). Let 
$$w_1=u(P_0), \, w_2=\displaystyle \frac{\overrightarrow{OA}}{|\overrightarrow{OA}|} \quad \textrm{and} \quad  
w_3=\alpha u(P_0)+ \beta u(P_1)$$ where 
\be
\label{eq:alpha_beta}
\alpha=\frac{-u(P_0) \cdot u(P_1)}{\sqrt{1-\big(u(P_0) \cdot u(P_1)\big)^2 }} \quad \textrm{and} \quad
 \beta=\frac{1}{\sqrt{1-\big(u(P_0) \cdot u(P_1)\big)^2 }}>0.
 \ee
The choice of $\alpha$ and $\beta$ is done in order to insure that $w_1 \cdot w_3=0$ and $|w_3|^2=1$ which finally yields
the orthonormal basis $w_1$, $w_2$ and $w_3$.
Note now that the vectors $u(P_0)$ and $u(P_1)$ have the following components in the basis $(w_1,w_2,w_3)$: 
$$
 u(P_0)=(1,0,0) \quad \text{ and } \quad u(P_1)=(-\frac{\alpha}{\beta},0,\frac{1}{\beta}).
$$
We want to find the writing of $\overrightarrow{P_0P_1}$ in that basis, too. For that, we have
$$
 \overrightarrow{P_0P_1}=  \overrightarrow{P_0A} + \overrightarrow{AB} +\overrightarrow{BP_1} 
$$ 
which implies the existence of three real numbers $\lambda, \tilde{\lambda}, \hat{\lambda}\in \RR$ with $\tilde{\lambda}\neq 0$ such that 
\begin{align*}
  \overrightarrow{P_0P_1} = & \lambda w_1+ \tilde{\lambda} w_2 + \hat{\lambda} u(P_1) \\
                = & \lambda w_1+ \tilde{\lambda} w_2 + \hat{\lambda} \left(\frac{1}{\beta}w_3-\frac{\alpha}{\beta}w_1 \right).
\end{align*}
Thus, $\overrightarrow{P_0P_1}$ has the following components in the basis $(w_1,w_2,w_3)$:
\[
 \overrightarrow{P_0P_1}=(\lambda-\frac{\alpha}{\beta}\hat{\lambda},\tilde{\lambda},\frac{\hat{\lambda}}{\beta}).
\]
Set now the following vector $\xi:=(1,s,-\beta)\neq 0$ written in our basis where $s:= \frac{\hat{\lambda}(\alpha+\beta)}{\beta\tilde{\lambda}}-\frac{\lambda}{\tilde{\lambda}}$. Then $[P_0P_1]\subset \Omega$ (since $\Omega$ is a ball) and
\begin{align*}
 & \overrightarrow{P_0P_1} \cdot \xi =0, \, \, \textrm{ i.e., } \, \, \xi\in P_0P_1^\perp\\
 & u(P_0) \cdot \xi = 1>0, \quad u(P_1) \cdot \xi = u(P_0) u(P_1)-1<0,
\end{align*}
which contradicts Proposition \ref{pro:preserv_prod}. Thus, $D_0$ and $D_1$ are coplanar.

\medskip

\begin{figure}[htbp]
\center
\includegraphics[scale=0.5,
width=0.5\textwidth]{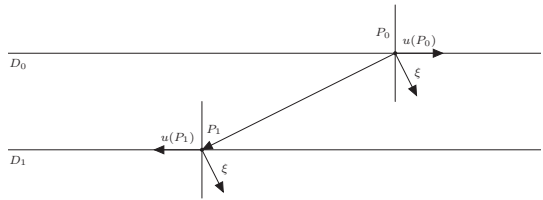} \caption{Two parallel lines $D_0$ and $D_1$}
\label{coplan} 
\end{figure}

\nd {\it Step 2: We show that $D_0$ and $D_1$ must intersect ($D_0$ might coincide with $D_1$)}. Assume by contradiction that $D_0$ and $D_1$ are parallel and $D_0\neq D_1$. It means that $u(P_0)=-u(P_1)$ 
(because of our choice $u(P_0)\neq u(P_1)$). Set $(w_1, w_2)$ be an orthonormal basis in the two-dimensional plane $\Pi$ determined
by $D_0$ and $D_1$ with $w_1=u(P_0)$. In the basis $(w_1, w_2)$, we write $\overrightarrow{P_0P_1}=(\lambda, \tilde{\lambda})$ where
$\tilde{\lambda}\neq 0$ (since $D_0\neq D_1$) and set $\xi=(-\tilde{\lambda}, \lambda)$ be an orthogonal vector to 
$\overrightarrow{P_0P_1}$ in $\Pi$ (see Figure \ref{coplan}). Then one checks that $u(P_0) \cdot \xi=-\tilde{\lambda}$ and $u(P_1) \cdot \xi=\tilde{\lambda}$ 
have different signs which again contradicts Proposition \ref{pro:preserv_prod}.

\medskip

\nd {\it Step 3: There exists a point $O\in D_0$ with $O\neq P_0, P_1$ and a sign $\gamma\in \{\pm 1\}$ such that $$u(P_i)=\gamma \frac{\overrightarrow{OP_i}}{|\overrightarrow{OP_i}|}, \quad i=0,1.$$} 
If $D_0= D_1$, then $u(P_0)=-u(P_1)$ so any point $O\in D_0$ located between $P_0$ and $P_1$ leads to the conclusion. 
Otherwise, $D_0\neq D_1$ and we denote $\{O\}=D_0\cap D_1$. First, we prove that $O\neq P_0, P_1$. Assume by contradiction that $O=P_0\in D_0\cap D_1$. Then by Proposition \ref{pro:trace} we know that $P_0$ and $P_1$ are Lebesgue points of the trace $\tilde{u}$ of $u$ on the segment $D_1\cap \Omega$ (directed by $u(P_1)$) with $\tilde{u}(P_0)=u(P_0)$ and $\tilde{u}(P_1)=u(P_1)$ so that by Lemma \ref{lem:const_droite},
we should have $u(P_0)$ is parallel with $u(P_1)$ which is a contradiction with $D_0\neq D_1$. So, $O\neq P_0, P_1$. Next, note that for any orthogonal vector $\xi$ to $\overrightarrow{P_0P_1}$
in the plane determined by $D_0$ and $D_1$, we have by Proposition \ref{pro:preserv_prod} that $u(P_0) \cdot \xi$ and $u(P_1) \cdot \xi$ have the same sign, i.e., 
\be
\label{eq:14_sign}
\big(u(P_0) \cdot \xi\big)\cdot \big(u(P_1) \cdot \xi\big)\geq 0.\ee Write now $$\overrightarrow{OP_0}=\lambda u(P_0)\quad \textrm{ and } \quad \overrightarrow{OP_1}=\tilde\lambda u(P_1)$$ with $\lambda, \tilde{\lambda}$ nonzero real numbers. The conclusion of Step 3 is equivalent with proving that $\lambda$ and $\tilde{\lambda}$ have the same sign. For that, 
as at Step 1, we choose the orthonormal basis
$w_1=u(P_0)$ and $w_2=\alpha u(P_0) +\beta u(P_1)$ with $\alpha\in \RR$ and $\beta>0$ given at \eqref{eq:alpha_beta} (recall that $|u(P_0) \cdot u(P_1)|<1$ because the assumption $D_0\neq D_1$). Since $\overrightarrow{P_0P_1}=\overrightarrow{OP_1}-\overrightarrow{OP_0}=\tilde{\lambda} u(P_1)-\lambda u(P_0)$, we write in the basis $(w_1, w_2)$:
\[
 u(P_0)=(1,0), \quad u(P_1)=(-\frac{\alpha}{\beta},\frac{1}{\beta}), \quad \overrightarrow{P_0P_1}=(-\lambda-\frac{\alpha}{\beta}\tilde{\lambda},\frac{\tilde{\lambda}}{\beta}).
\]
Then for the orthogonal vector $\xi=(\tilde{\lambda}, \lambda\beta+\alpha\tilde{\lambda})\neq 0$ to $\overrightarrow{P_0P_1}$ we have by \eqref{eq:14_sign}:
$$0\leq \big(u(P_0) \cdot \xi\big)\cdot \big(u(P_1) \cdot \xi\big)= \tilde{\lambda}\cdot {\lambda}.
$$

\medskip

\nd {\it Step 4: Conclusion}. 
For every Lebesgue point $P\in Leb\cap \Omega$ of $u$, we consider the line $D$ passing through $P$ and directed by $u(P)$. Call $\mathcal{D}$ the set of these lines. Obviously, $\mathcal D$ covers $\h^N$-almost all the ball $\Omega$ (since $\h^{N}(\Omega\setminus Leb)=0$), in particular, $\cal D$ is not planar. By Step 1, we know that every two lines in $\mathcal D$ are coplanar. Then Proposition \ref{prop:geometry} (whose proof is presented below) implies that either all these lines are parallel, or they pass through the same point $O$. Since $u$ is nonconstant, we deduce by Step 2 that only the last situation holds true. By Step 3, we conclude that 
$u=\gamma \us(\cdot-O)$ a.e. in $\Omega$.

\medskip

\nd {\bf Case II: $\Omega$ is a connected open set}. By Case I, we know that in every open ball $B\subset \Omega$ around a Lebesgue point of $u$, the vector field $u$ is either a vortex type vector field in $B$, or $u$ is constant in $B$. Since $u$ is nonconstant in $\Omega$, there exists a Lebesgue point $P_0$ of $u$ and a ball $B_0\subset \Omega$ around $P_0$ such that $u$ is a vortex type vector field in $B_0$, say for simplicity $u=\us$. Let $P\neq P_0$ be any other Lebesgue point of $u$. Since $\Omega$ is path-connected, there exists a path $\Gamma\subset \Omega$ from $P_0$ to $P$. Then we can cover the path $\Gamma$ by a finite number of open balls $\{B_j\}_{0\leq j\leq n}$ such that $P\in B_n$, $B_j\cap B_{j+1}\neq \emptyset$ for $0\leq j\leq n-1$ and $u$ is either constant, or a vortex type vector field in any $B_j$. Since $u=\us$ in $B_0$ and $B_0\cap B_1$ is a nonempty open set, the analysis in Case I yields $u=\us$ in $B_1$ and by induction, $u=\us$ in $B_n$ which is a neighborhood of $P$. This concludes our proof.
\qed

\medskip

Let us now present the proof of the geometric result in Proposition \ref{prop:geometry} which is independent of the previous results:

\proof{ of Proposition \ref{prop:geometry}}
Assume that there are two lines $D_0, D_1\in \mathcal{D}$ that are not colinear. Since $D_0$ and $D_1$ are coplanar, they intersect at a point $P$.
Call $\Pi$ the plane determined by $D_0$ and $D_1$. We show that all the lines in $\mathcal D$ pass through $P$. Let $D_2\in {\mathcal D}$ be 
any line not included in $\Pi$ (such line exists because $\mathcal D$ is not planar). We know that $D_2$ is coplanar with $D_0$ and $D_1$, respectively. Then $D_2$ cannot be parallel with $D_0$ (otherwise, $D_2\parallel D_0$ and $D_2\cap D_1\neq \emptyset$ imply that $D_2\subset \Pi$ which is a contradiction with our assumption). Similarly, $D_2$ cannot be parallel with $D_1$. Therefore, $D_2$ intersects both $D_0$ and $D_1$. Since $D_2$ is not included in $\Pi$, the intersection points coincide with $P$. Let now $D_3\in \mathcal{D}$ be any line included in $\Pi$ (different than $D_0$ and $D_1$). Then $D_3$ is not included in the plane determined by $D_1$ and $D_2$. The previous argument leads again to $P\in D_3$ which concludes our proof.
\qed

\section{Vortex-line type vector fields}\label{sec:weak_form}
We will prove the characterization of the weaken kinetic formulation \eqref{eq:kineNW} in Theorem \ref{thm:vortex-line}. This result is in the spirit of 
Corollary \ref{cor:vortex} and leads to vector fields that have vortex-line singularities.

\proof{ of Theorem \ref{thm:vortex-line}}
For $x \in \R^N$, recall the notation $x=(x', x_N)$ with $x'=(x_1,\ldots,x_{N-1}) \in \R^{N-1}$. As the result is local in the set $\{u_N\neq \pm 1\}$, we will assume that $\om=B'\times (-1,1)$ is included in that set where $B'$ is the unit ball in $\R^{N-1}$. Let $\xi' \in \Ss^{N-2}$ and $\xi=(\xi',0) \in \mathcal{E}$. 
Since $e_N \in \xi^{\perp}$, we deduce by \eqref{eq:kineNW}:
\begin{equation}\label{eq:constante_e_n}
 e_N \cdot \nabla_x \chi(.,\xi)=\partial_N \chi(., \xi)=0 \text{ in $\mathcal{D}'(\omega)$}.
\end{equation}
We know that the point $(x',t)$ is a Lebesgue point of $\chi(., \xi)$ for $\h^{N-1}$-a.e. $x'\in B'$ and $\h^1$-a.e. $t\in (-1,1)$ and the convolution argument in the proof of Proposition \ref{pro:preserv_prod} yields
$$\chi(x,\xi)=\chi(x+te_N, \xi) \quad \textrm{for $\h^{N}$-a.e. $x\in \om$ and $\h^1$-a.e. $t$}.$$
Then one can define the measurable function $\tilde \chi(\cdot, \xi'):B'\to \{0,1\}$ by 
$$\tilde{\chi}(x',\xi'):=\chi(x,\xi)=\mathds{1}_{\{x\in \om\, :\, u'(x) \cdot \xi'>0\}} \quad \textrm{for $\h^{N}$-a.e. $x=(x',t)\in \om$}.$$
Set
\[
 \tilde u(x')=\frac1 {V_{N-2}}\int_{\Ss^{N-2}} \xi' \tilde{\chi}(x',\xi') \diff \mathcal{H}^{N-2}(\xi'), \quad x'\in B'.
\]
Thanks to \eqref{eq:averN}, $$\tilde u(x')=\frac{u'(x)}{|u'(x)|}\quad \textrm{ for $\mathcal{H}^N$-a.e. $x=(x',t)\in \om \subset \{ |u'|>0 \}$}.$$ 
In particular, $\tilde{\chi}(x',\xi')=\mathds{1}_{\{x'\in B'\, :\, \tilde u(x') \cdot \xi'>0\}}$ in $B'$ for every $\xi'\in \Ss^{N-2}$. Therefore, we deduce by \eqref{eq:kineNW} that 
$\tilde{u} : B' \to \Ss^{N-2}$ satisfies:
\[
 \forall \xi' \in \Ss^{N-2}, \; \forall v' \in (\xi')^{\perp}, \quad v' \cdot \nabla'_{x'} \tilde{\chi}(x',\xi')=0 \text{ in } B'.
\]
where $\nabla'_{x'}=(\partial_1, \dots, \partial_{N-1})$. As $N-1\geq 3$, Theorem \ref{thm:main} yields either $\tilde u(x')=w'$ for almost every $x'\in B'$ where $w' \in \Ss^{N-2}$ is a constant vector, or 
$\tilde{u}(x')=\gamma \frac{x'-P'}{|x'-P'|}$ for almost every $x'\in B'$ where $\gamma\in \{\pm 1\}$ and $P' \in \R^{N-1}$ is some point.
This means that for a.e. $x \in \omega$,
\[
\text{either} \quad u'(x)=|u'(x)|w' \quad \text{or} \quad u'(x)=\gamma |u'(x)|\frac{x'-P'}{|x'-P'|}.
\]

\medskip

\nd {\it Case 1.} Let $u'(x)=|u'(x)|w'$ for  a.e. $x \in \omega$.
By \eqref{differ_deriv}, we have for $k \in \{1,\ldots,N-1\}$,
\begin{align} \label{eq:permut}
\partial_k u_N =  \partial_N u_k  = w_k \partial_N (|u'|) \quad \textrm{ in } \om
\end{align}
which yields for all $k,j \in \{1,\ldots,N-1\}$:
\begin{align*} w_j\partial_k u_N=w_k \partial_{j} u_N  \quad \textrm{ in } \om.
\end{align*}
Therefore, $u_N(x)=g(\alpha, x_N)$ in $\om$ for some two-dimensional function $g$ with the new variable $\alpha:=\alpha(x)=x'\cdot w'$.
Moreover, by \eqref{eq:permut}, the function $g$ satisfies the following: since $w_k\neq 0$ for some $k\in \{1,\ldots,N-1\}$
(because $w\in \Ss^{N-1}$) then the equation $|u'|^2+u_N^2=1$ a.e. in $\om$ implies
$$w_k \partial_\alpha g=\partial_k  u_N\stackrel{\eqref{eq:permut}}{=}w_k \partial_N (|u'|)=w_k \partial_N (\sqrt{1-g^2}).$$
The Poincar\'e lemma yields the existence of a stream function $\psi(\alpha, x_N)$ such that $g=\partial_N \psi$ and $\sqrt{1-g^2}=\partial_\alpha \psi$ so that
$u(x)=\nabla_x[\psi(\alpha, x_N)]$ and therefore, $\psi$ satisfies
the two-dimensional eikonal equation:
$$(\partial_\alpha \psi)^2+(\partial_N \psi)^2=1.$$

\medskip

\nd {\it Case 2.} Let $u'(x)=\gamma |u'(x)|\frac{x'-P'}{|x'-P'|}$ for  a.e. $x \in \omega$. As above we have for $k \in \{1,\ldots,N-1\}$:
\begin{align}\label{eq:permut2}
\partial_k u_N  =  \partial_N u_k  = \gamma \frac{x_k-P_k}{|x'-P'|} \partial_N (|u'|)   \quad \textrm{ in } \om
\end{align}
and we deduce that for all $k,j \in \{1,\ldots,N-1\}$:
\begin{align*}
 (x_j-P_{j}) \partial_k u_N=(x_k-P_k) \partial_{j} u_N \quad \textrm{ in } \om.
\end{align*}
Therefore, $u_N(x)=g(\alpha, x_N)$ in $\om$ for some two-dimensional function $g$ with the new variable $\alpha:=\alpha(x)=|x'|$.
By \eqref{eq:permut2}, we conclude as above that there exists a stream function $\psi$ solving the eikonal equation in the variables $(\alpha, x_N)$ such that
$$u(x)=\nabla_x [\psi(\alpha, x_N)].$$ \qed

\bigskip

\nd

\textbf{Acknowledgments}: RI acknowledges partial support by the ANR project ANR-14-CE25-0009-01.

\bibliographystyle{plain}
\bibliography{bib}
\end{document}